\documentclass{article}

\usepackage{amsmath,amsfonts,amsthm,amssymb,amscd,color,xcolor,mathrsfs,verbatim,microtype}
\usepackage{graphicx,eurosym}
\usepackage{hyperref}

\usepackage[applemac]{inputenc}

\usepackage[cyr]{aeguill}

\colorlet{darkblue}{blue!50!black}

\hypersetup{
    colorlinks,%
    citecolor=blue,%
    filecolor=red,%
    linkcolor=darkblue,%
    urlcolor=blue,%
    pdfnewwindow=true,%
    pdfstartview={FitH}
}

\usepackage{graphicx,amscd,mathrsfs,wrapfig,mathrsfs,lipsum}
\usepackage{eufrak}
\usepackage{float}
\usepackage{tikz}
\usepackage{multicol}
\usepackage{caption}
\usetikzlibrary{arrows}
\usepackage{capt-of}

\colorlet{darkblue}{blue!50!black}

\newcommand{\p}{\partial}
\newcommand{\e}{\varepsilon}

\newcommand{\ppP}{{\mathsf P}}

\newcommand{\R}{{\mathbb R}}
\newcommand{\Z}{{\mathbb Z}}
\newcommand{\IP}{{\mathbb P}}
\newcommand{\pP}{{\mathbb P}}
\newcommand{\I}{{\mathbb I}}

\newcommand{\T}{{\mathbb T}}

\newcommand{\N}{{\mathbb N}}
\newcommand{\la}{\lambda}

\newcommand{\La}{\Lambda}
\newcommand{\ty}{\infty}
\newcommand{\te}{\theta}

\newcommand{\de}{\delta}

\newcommand{\bbar}{\boldsymbol{|}}

\newcommand{\Ker}{\mathop{\rm Ker}\nolimits}

\newcommand{\z}{\mathfrak{z}}

\newcommand{\h}{{\mathfrak h}}

\newcommand{\aA}{{\cal A}}
\newcommand{\BB}{{\cal B}}

\newcommand{\EE}{{\cal E}}
\newcommand{\FF}{{\cal F}}
\newcommand{\GG}{{\cal G}}
\newcommand{\HH}{{\cal H}}

\newcommand{\KK}{{\cal K}}

\newcommand{\OO}{{\cal O}}
\newcommand{\PP}{{\cal P}}

\newcommand{\bB}{{\mathfrak{b}}}

\newcommand{\RR}{S}

\newcommand{\VV}{{\cal V}}

\newcommand{\XX}{{\cal X}}

\newcommand{\lag}{\langle}
\newcommand{\rag}{\rangle}

\newcommand{\dd}{{\textup d}}

\newcommand{\PPPP}{{\mathfrak P}}

\newcommand{\supp}{\mathop{\rm supp}\nolimits}
\newcommand{\diver}{\mathop{\rm div}\nolimits}

\theoremstyle{plain}

\newtheorem*{mta}{Theorem~\hypertarget{A}{A}}
\newtheorem*{mtb}{Theorem~\hypertarget{B}{B}}
\newtheorem*{mtc}{Theorem~\hypertarget{C}{C}}

\newtheorem*{lemma*}{Lemma}
\newtheorem{theorem}{Theorem}[section]
\newtheorem{lemma}[theorem]{Lemma}
\newtheorem{proposition}[theorem]{Proposition}
\newtheorem{corollary}[theorem]{Corollary}
\theoremstyle{definition}
\newtheorem{definition}[theorem]{Definition}

\theoremstyle{remark}

\newtheorem{remark}[theorem]{Remark}

\numberwithin{equation}{section}

\begin{document}
\author{Pierre-Marie Boulvard\footnote{Institut de Math\'emathiques de Jussieu--Paris Rive Gauche, Universit\'e Paris Diderot, UMR 7586, Sorbonne Paris Cit\'e, F-75013, Paris, France; e-mail: \href{mailto:pierre-marie.boulvard@imj-prg.fr}{Pierre-Marie.Boulvard@imj-prg.fr}} \and
Peng Gao\footnote{School of Mathematics and Statistics, Center for Mathematics and Interdisciplinary Sciences, Northeast Normal University, Changchun 130024, China;
e-mail: \href{mailto:gaop428@nenu.edu.cn}{gaop428@nenu.edu.cn}}
 \and Vahagn~Nersesyan\footnote{NYU-ECNU Institute of Mathematical Sciences, NYU Shanghai, 3663 Zhongshan Road North, Shanghai, 200062, China, e-mail: \href{mailto:vahagn.nersesyan@nyu.edu}{Vahagn.Nersesyan@nyu.edu}}\, \footnote{Universit\'e Paris-Saclay, UVSQ, CNRS, Laboratoire de Math\'ematiques de Versailles, 78000, Versailles, France}}

  \title{Controllability and ergodicity of
3D primitive equations driven by a    finite-dimensional~force}
\date{\today}

\maketitle

\begin{abstract}

We study the problems of    controllability and ergodicity of the system of 3D primitive equations modeling large-scale oceanic and atmospheric motions.~The system is driven by  an additive     force acting   only   on a finite number of  Fourier   modes in the temperature equation.~We~first show  that  the     velocity and   temperature components of the equations    can be simultaneously approximately controlled to     arbitrary position in the phase space. The~proof is based on     Agrachev--Sarychev type geometric control~approach.

Next, we   study   the  controllability of the    linearization  of~primitive equations  around a non-stationary    trajectory of the randomly forced system. Assuming that the probability~law of the forcing is decomposable and observable, we prove  almost sure approximate controllability     by using the same   Fourier modes as in the nonlinear setting.~Finally,   combining  the      controllability  of the  linearized system with   a criterion from~\cite{KNS-2018}, we~establish exponential mixing for  the  nonlinear    primitive~equations with a  random~force.


\end{abstract}

\newpage
 \tableofcontents

\setcounter{section}{-1}

\section{Introduction}
\label{S:0}

The system of {\it 3D primitive equations}  (PEs) of meteorology and oceanology is an important model   of  geophysical fluid dynamics. Today, most numerical weather prediction and climate simulation models are based on them. This system is derived, using hydrostatic approximation, from the 3D~Navier--Stokes equations with Coriolis force coupled with the thermodynamic equation  (see the book by Zeitlin~\cite{zeitlin2018geophysical}). The mathematical study of these equations has attracted a lot of attention in the last two decades.~Following the   framework introduced 
by Lions, Temam, and  Wang~\cite{LTS-92b, LTS-92a}, we consider in this paper  the PEs     written in the form:
\begin{gather}
     \p_t v + L_1v+ \lag v, \nabla\rag v + w\,
\p_z v   +   f   v^\bot   +
\nabla p = h_1,  \label{0.1}  \\
 \p_z p  + \te   =0,    \label{0.2}  \\
  \diver  v +\p_z w =0,   \label{0.3} \\
 \p_t \te   + L_2 \te  + \lag v,  \nabla\rag \te + w\, \p_z\te 
    =  h_2+\eta. \label{0.4}
\end{gather}
The unknowns are the 3D velocity field of the fluid $(v_1,v_2,w)$, where\footnote{We denote $v^\bot=(-v_2,v_1)$.} $v=(v_1,v_2)$  and~$w$ are  the horizontal     and vertical  velocity components, the temperature~$\theta$, and  the pressure~$p$.~The number~$f$ is the Coriolis rotation frequency,    the functions~$h_1$ and $h_2$  are given    source terms, and~$\eta$   is an external perturbation---a~{\it control}   or  a~{\it random noise}.  The~operators 
\begin{align*}
L_1&=-\nu_1\Delta-\mu_1 \p_{zz} ,\\
L_2&=-\nu_2\Delta-\mu_2 \p_{zz}	
\end{align*}are the viscosity and   heat diffusions, where the   numbers $\nu_1,\mu_1>0$ are the horizontal and vertical viscosities, while    $\nu_2,\mu_2>0$ are the horizontal and vertical heat diffusivity coefficients.~We denote by $ \Delta, \nabla$,   $\diver$ the 2D (horizontal) Laplacian, gradient,   divergence operators: 
$$
\Delta  = \p_{xx}  + \p_{yy}  , \quad \nabla  = \left(\p_x, \p_y\right), \quad  \diver  = \lag \nabla , \cdot \rag,
 $$
and $\lag v, \nabla\rag = v_1 \p_x +v_2\p_y$.

     The space variable $  (x,y,z)$ is assumed to belong  to the torus $\T^3 =\R^3/2\pi\Z^3$, i.e., all the  above functions are  $2\pi$-periodic in  $x,y,$ and~$z$. Furthermore,  we   assume that the functions $v, p,  h_1$ are even and the functions $w, \theta, h_2, \eta$ are odd   in~$z$. As a consequence,~$w, \theta, h_2,  \eta$~vanish at~$z=0$.
     
     The unknown functions in system~\eqref{0.1}-\eqref{0.4} can be divided into two types: the prognostic unknowns $v$ and $\theta$, which are determined through  an initial boundary value problem, and diagnostic ones $w$ and $p$, which can be expressed as functions of $v$ and $\theta$.~Indeed, from the conservation of   mass equation \eqref{0.3} and the boundary condition $w|_{{z=0}}=0$  it follows that   
    \begin{equation}\label{0.5}
     w(t,x,y,z)= -\int_0^z \diver  v(t,x,y,\z)\,\dd \z, 
        \end{equation}
and from the hydrostatic balance \eqref{0.2} that
\begin{equation}\label{0.6}
p(t,x,y,z)=p_s(t,x,y)-\int_0^z \te(t,x,y,\z)\,\dd \z.
\end{equation}Using   equalities~\eqref{0.5} and~\eqref{0.6},     the following equivalent formulation is obtained for the~PEs:
\begin{gather}
  \p_t v + L_1v+ \lag v, \nabla\rag v - \int_0^z \diver  v(t,x,y,\z)\,\dd \z \,
\p_z v  +   f   v^\bot \nonumber   \\+
\nabla p_s(t,x,y)-\int_0^z \nabla\te(t,x,y,\z)\,\dd \z = h_1,  \label{0.7}   \\
   \p_t \te   + L_2 \te  + \lag v,  \nabla\rag \te -  \int_0^z \diver v(t,x,y,\z)\,\dd \z\,
\p_z \te     =  h_2+\eta. \label{0.8}
\end{gather}The     well-posedness of these equations     has been studied by many authors.  The existence of {\it weak  solutions} is known from the works of Lions, Temam, and  Wang~\cite{LTS-92b, LTS-92a}, but the  uniqueness   is still an open problem. In~this paper, we  deal with   {\it strong solutions}   whose  global  existence and uniqueness is established by Cao and  Titi~\cite{CT-2007} in the case of   Neumann boundary conditions; see also the paper by Kobelkov~\cite{Kob-07} for a different proof.  In~the case of periodic boundary conditions, the global existence of strong solutions  is considered      by Petcu~\cite{P-2006}   and~in~the case of    Dirichlet boundary  conditions,   by Kukavica and   Ziane~\cite{KZ-2007}. The~existence of a global attractor is obtained by Ju~\cite{J-2007} and Chueshov~\cite{Chu-14}. We refer the reader   to the    reviews \cite{TZ-2004,PTZ-2009} for more details  and references.

  In the periodic setting, the  PEs~\eqref{0.7},~\eqref{0.8}   are considered in the     function spaces~$H$ and~$V$       recalled in Section~\ref{S:1}.  
  To formulate the first main result of this paper, we assume that 
  the couple of source terms~$(h_1,h_2)$ is a smooth element of~$H$, and  $\eta$ is a control taking           values in the space $\HH= \text{span} \{\phi_i: i=1,\ldots,10\}$, where $\phi_i$ are   the following eigenfunctions of the  heat   diffusion   operator~$L_2$:   
  $$
\cos j x\sin z,~~ \sin jx\sin z,~~\cos j y\sin z,~~ \sin j y \sin z,~~  \sin j z, \quad j=1,2.
$$   
\begin{mta}
 Problem \eqref{0.7},~\eqref{0.8} is approximately controllable by $\HH$-valued controls.~More precisely, for any~$\e>0$, any time $T>0$,   any initial condition~$(v_0,\theta_0)\in V$, and any target $(v_1,\theta_1)\in H$, there is a control $\eta\in L^\ty([0,T],\HH)$ such that the unique strong solution~$(v,\te)$  of problem~\eqref{0.7},~\eqref{0.8}   satisfies    
\begin{gather}
	(v(0), \te(0))=(v_0,\theta_0),\label{0.9}
	\\ \|(v(T), \te(T))-(v_1,\theta_1)\|_{L^2(\T^3,\R^3)}<\e.\label{0.10}
\end{gather}
 \end{mta}  
 Note that the space $\HH$ of admissible values for the control $\eta$ is {\it independent} of the physical  parameters  $h_1,h_2, f, \nu_i,\mu_i$, $i=1,2$.~A more general version of this result is formulated in Theorem~\ref{T:2.3}, where  a saturation property is specified  that ensures the approximate controllability of the system. We also show that if  some   controlled Fourier modes are added in   the velocity   equation~\eqref{0.7},  then    approximate  controllability   holds with  respect to the stronger norm of the space~$H^1(\T^3,\R^3)$.

 Approximate controllability of PDEs by additive finite-dimensional forces has been studied by many authors in the recent years.~The first results are obtained by Agrachev and Sarychev~\cite{AS-2005, AS-2006}, who considered     the Navier--Stokes~(NS)  and Euler systems on the 2D torus (see also the review \cite{AS-2008}).~Their approach has been    generalized  by    Shirikyan~\cite{shirikyan-cmp2006, shirikyan-aihp2007} to the case of    the   3D NS     system; see also the papers~\cite{Shi-2013,Shir-2018} by Shirikyan, where     the Burgers equation is considered   on the real line     and      on a bounded  interval with   Dirichlet boundary conditions.  
   Rodrigues and Phan~\cite{SSRodrig-06, RD-2018} established  approximate     controllability    of the NS system  on  2D and 3D rectangles with   Lions boundary conditions.~In the periodic setting, Nersisyan~\cite{Hayk-2010,   Hayk-2011} considered    3D Euler systems for perfect compressible and incompressible     fluids,     Sarychev~\cite{Sar-2012} studied    the  2D   cubic Schr\"odinger equation, and  Nersesyan~\cite{VN-2015} considered the     Lagrangian trajectories of the~3D~NS~system.

The proof of Theorem~\hyperlink{A}{A}   is based on a technique of applying    large controls on short time intervals.   Previously, such ideas have been used  mainly   in the study of  finite-dimensional control systems; e.g., see the works of Jurdjevic and  Kupka~\cite{JK-1985, MR1425878} and the references therein. 
 Infinite-dimensional   extensions of this technique  appear     in the   above-cited papers of Agrachev and Sarychev. More recently, this approach has  been used  in the paper  of Glatt-Holtz, Herzog, and Mattingly~\cite{GHM-2018}, where, in particular,  a   1D
parabolic PDE  is considered with polynomial nonlinearity of odd degree, and  in the paper of Nersesyan~\cite{nersesyan-2018}, where the nonlinearity  is a smooth function that grows   polynomially  without any restriction on the degree   and on the space dimension.

The main difficulty of the problem considered in this paper comes from the {\it highly degenerate} nature of the control system. The form of the saturation property and the argument for its verification   are    more complicated than in the  previously studied situations. When the control acts directly only on the temperature equation, we are able to check the saturation       with respect to the $L^2$-norm. The latter  is known   to be poorly adapted     for the stability properties of the 3D PEs and is a source of many difficulties  in different parts of the proof. %This loss of stability is compensated by a use of   more advanced  and finer perturbative analysis.  

  % + Inch vor tex nshel vor bnakan a henc karavarel jermutyun@+hxumner

%+ Normal bacatrel inchqan enq \cite{nersesyan-2018}-ic ogtvum: }

To formulate our second   result, let us assume that~$\eta$ is a 
{\it Haar coloured noise} taking values in the same space $\HH$ as above.   This means that $\eta$ has  the form  
\begin{equation} \label{0.11}
	\eta(t)=  \sum_{i=1}^{10}   \eta^i(t)\phi_i,
\end{equation}
where~$\{\eta^i\}$ are   independent copies of a random process $\tilde\eta$  defined by
\begin{equation} \label{0.12}
	\tilde\eta(t)=\sum_{k=0}^\infty \xi_k \h_0(t-k)
	+\sum_{j=1}^\infty j^{-q}\sum_{l=0}^\infty \xi_{jl}\h_{jl}(t).
\end{equation}
Here $q>1$,    $\{\xi_k,\xi_{jl}\}$  are independent identically distributed (i.i.d.)~scalar random variables with Lipschitz-continuous density~$\rho$, and  $\{\mathfrak{h}_0,\h_{jl}\}$ is the Haar system defined by  
\begin{align*}
	h_0(t)&=\left\{
	\begin{array}{cl}
		1 & \mbox{for $0\le t<1$},\\[2pt]
		0 & \mbox{for $t<0$ or $t\ge1$},
	\end{array}
\right.
\\
	h_{jl}(t)&=\left\{
	\begin{array}{cl}
		0 & \mbox{for $t<l2^{-j}$ or $t\ge (l+1)2^{-j}$},\\[2pt]
		1 & \mbox{for $l2^{-j}\le t<\bigl(l+\tfrac12\bigr)2^{-j}$},\\[2pt]
		-1 & \mbox{for $\bigl(l+\tfrac12\bigr)2^{-j}\le t<(l+1)2^{-j}$}
	\end{array}
\right.
\end{align*}
with $j\ge1$ and $0\le l\le 2^j-1$. Note that, for any   $k\ge1$, the functions
$$
\bigl\{h_0(\cdot-k+1),\,h_{jl},\, j\ge1,\,(k-1)2^j\le l\le k2^j-1\bigr\} 
$$
are supported on the interval~$[k-1,k]$ and form an orthonormal basis in the space~$L^2([k-1,k])$.

  Let     $(\tilde v,\tilde \te)$ be a trajectory  of problem~\eqref{0.7},~\eqref{0.8} with process  $\eta$ defined by~\eqref{0.11}. 
  \begin{mtb}Under the above conditions,
 	the   linearization of problem \eqref{0.7},~\eqref{0.8} around the trajectory $(\tilde v,\tilde \te)$  (see system \eqref{1.11})  is  almost surely approximately controllable by $\HH$-valued controls.
 \end{mtb}
 See Theorem~\ref{T:3.3} for a   more precise formulation of the result. It is    proved by showing that the kernel of the {\it random  Gramian operator} is almost surely non-trivial.~The latter is derived from   the {\it observability~property} of the Haar~noise. Let~us  emphasize that   on a non-empty, but zero-probability event (depending on the reference trajectory~$(\tilde v,\tilde \te)$),   the  linearized problem   is non-controllable.~Indeed,  
  assume that the source terms~$h_1$ and~$h_2$, the noise~$\eta$, as well as  the   trajectory~$(\tilde v,\tilde \te)$,  are identically zero.~Then the space~$\HH$ is invariant for the    linearized problem, so the attainable  set from the origin  cannot be dense in $H$.

 Controllability   properties  of    nonlinear and  linearized    equations have     applications to the study of   randomly perturbed problems. Indeed, it is well known that   approximate controllability     implies, for example,   irreducibility of the associated Markov process  when the support of the law of the   noise is the whole space~$L^\ty([0,T],\HH)$ (see Section~6.3 in~\cite{KS-book}   for more details).~As it is shown in the recent papers by Kuksin, Nersesyan, and Shirikyan~\cite{KNS-2018, KNS-2019}, 
    the controllability of the  linearized system  can be used      in the analysis of the ergodicity problem  when the system is perturbed by a bounded degenerate  noise.      In~these papers   the      NS system, complex Ginzburg--Landau equations, and  parabolic PDEs with polynomial nonlinearities are studied. 
    See also the papers~\cite{KZ-2018} and~\cite{vnersesyan-2019}   for some related  situations where the noise is     non-degenerate.~In our third result,
   we show that the approach of these papers can be extended to the more degenerate  case of PEs.     
  To formulate the result,  let~$((v_k,\te_k), \pP_{(v,\te)})$~be the Markov family obtained by restricting   the trajectories of system \eqref{0.7},~\eqref{0.8},~\eqref{0.11} to integer times. Recall that $\rho$ is the density of the random variables $\{\xi_k,\xi_{jl}\}$  in~\eqref{0.12}.

   \begin{mtc} In addition to the above conditions, assume that  $(h_1,h_2)=0$,  the support~of~the density $\rho$ is bounded, and     $\rho(0)>0$.~Then
   the family $((v_k,\te_k), \pP_{(v,\te)})$ has a
unique stationary measure  on $V$  which is 
  exponentially mixing
in the dual-Lipschitz metric.
  \end{mtc}

See Section~\ref{S:mix} for more details. In the case of  the 3D~primitive equations with spatially regular white noise,   existence of   stationary measure is established by Glatt-Holtz,   Kukavica,   Vicol, and  Ziane~\cite{GKVZ-2014}. As far as we know,   uniqueness of stationary measure for that situation is still an open~problem due to    rather weak tail estimates for  
solutions.  The boundedness of the noise allows to reduce the study of the system to a compact phase space.~This naturally eliminates the problems coming from  the  tail estimates. On the other hand, bounded noises are well-accepted and commonly used in the physics literature (e.g., see~\cite{Ono13} and the references therein).
 In the case of   {\it non-degenerate} bounded kick force,   uniqueness   and   exponential mixing are proved by Chueshov~\cite{Chu-14}.  Let us also recall    some previous results considering  the problem of ergodicity for    PDEs driven by  a   degenerate noise.
      Hairer and Mattingly~\cite{HM-2006, HM-2011} used   Malliavin calculus      to study the  ergodicity for the NS system with   a white noise which  is degenerate in the Fourier space.  F\"oldes, Glatt-Holtz, Richards, and
Thomann~\cite{FGRT-2015} considered a similar problem for the Boussinesq system. Using controllability methods, Shirikyan~\cite{shirikyan-asens2015, shirikyan-2018} studied the NS system with a noise that is   localized in the physical space (distributed   in a subdomain or on the boundary). For   more results and references,  we refer the reader to   the book~\cite{KS-book}.

 This paper is  organized as follows.~In Section~\ref{S:1}, we   recall the functional setting for the PEs and formulate    perturbative results with respect to the initial condition and   control. In~Sections~\ref{S:2} and~\ref{S:3}, we discuss the problems of controllability of   nonlinear and  linearized  PEs and prove  Theorems~\hyperlink{A}{A} and~\hyperlink{B}{B}. In~Section~\ref{S:mix}, we consider the randomly forced PEs and  prove Theorem~\hyperlink{C}{C}.
  Examples of saturating spaces are provided in Section~\ref{S:4}. Finally, in   Section~\ref{S:6}, we  establish a   perturbative result  formulated in       Section~\ref{S:1}.

  \subsubsection*{Acknowledgement}
 
The authors thank the reviewers for their valuable suggestions.~They also thank Sergei Kuksin, Madalina Petcu, and Armen Shirikyan for discussions.  
  PG would like to thank the China
Scholarship Council (No.~201806625036) for the financial support    and the CNRS and
IMJ, Universit\'e Paris Diderot-Paris~7 for the hospitality   during his visit from December 2018 to
November 2019.
  The~research of VN was supported by the   ANR through  the grant NONSTOPS   ANR-17-CE40-0006-02 and by the CNRS through  the PICS grant {\it Fluctuation theorems in stochastic~systems}.

\subsubsection*{Notation}

  Throughout this paper, we   use the following notation.

  \smallskip
\noindent
$\Z^d,~d\ge1$ is   the integer lattice in~$\R^d$, and  $\T^d$ is the     torus~$\R^d/2\pi\Z^d$.

\smallskip
\noindent
$L^p(\T^d,\R^n),~p\ge 1,~n\ge1$ and $H^k(\T^d,\R^n),~k\ge 0$ are the usual Lebesgue~and  Sobolev spaces of functions $g:\T^d\to \R^n$  endowed with the    norms~$\|\cdot\|_{L^p}$ and~$\|\cdot\|_k$,~respectively. If $p=2$, we write $\|\cdot\|$ instead of $\|\cdot\|_{L^2}$ and denote by~$\lag\cdot,\cdot\rag$ the corresponding scalar product. 
 If $p=+\ty$, we write $\|\cdot\|_\ty$ instead of~$\|\cdot\|_{L^\ty}$.

\smallskip
   \noindent
$C^\ty(\T^d,\R^n)$ is the space of  infinitely differentiable functions $g:\T^d\to \R^n$.

\smallskip

    \noindent
Let   $X$ be a Banach  space endowed with the
norm $\|\cdot\|_X$.

  \smallskip
\noindent $B_X(a,r)$ denotes  the closed ball of radius $r>0$ centred at $a\in X$.

\smallskip
\noindent
$\BB(X)$ is the Borel $\sigma$-algebra on~$X$, and
 $\PP(X)$ is the set of Borel probability measures on $X$.

\smallskip
\noindent  $L^p(J_T,X)$, $1\le p<\ty$         is  the
space of measurable functions $u: J_T\rightarrow X$ endowed with the~norm
$$
\|u\|_{L^p(J_T,X)}=\left(\int_0^T \|u(t)\|_X^p\,\dd t
\right)^{1/p}, \quad J_T=[0,T]. 
$$

  \noindent
$C(J_T,X)$ (resp. $L^\ty(J_T,X)$) is the space of continuous (resp. bounded measurable) functions $u:J_T\to X$   endowed with the~norm
$$
\|u\|_{C(J_T,X)} (\text{resp. } \|u\|_{L^\ty(J_T,X)})=\sup_{t\in J_T} \|u(t)\|_X.
$$

 \section{Preliminaries on primitive equations}\label{S:1}
 
  We consider the system of PEs  
 in the     spaces $H^k$, $k\ge0$  defined by
$$    
H^k = \text{closure of $\VV$ in $H^k(\mathbb{T}^{3},\R^3)$}
$$ and endowed with the Sobolev norms $\|\cdot\|_k$ (with $L^2$-norm $\|\cdot\|$ if $k=0$),
   where~$\VV=\VV_1\times \VV_2$ and $\VV_1$ and $\VV_2$ are the spaces given by
  \begin{align*}
\VV_1&=\left\{v\in   C^\ty(\mathbb{T}^{3},\R^2): v~\text{is even in}~z, \int_\T\diver v\,\dd z=0, \int_{\T^3} v \,\dd x\,\dd y\,\dd z=0 \right\}, \\
\VV_2&=\left\{\te\in  C^\ty(\mathbb{T}^{3},\R): \te~\text{is odd in}~z, \int_{\T^3} \te \,\dd x\,\dd y\,\dd z=0\right\}.\end{align*}
The condition $\int_\T \diver v\,\dd z =0$ in the definition of~$\VV_1$ comes from equality~\eqref{0.5} at~$z=2\pi$, the fact that $w$ is $2\pi$-periodic, and the boundary value~$w|_{z=0}=0$;   see~\cite{CT-2007,PTZ-2009} for more details. We will mainly consider the spaces\footnote{The subscripts $1$ and $2$ are used with
 $H,V, U$ to denote spaces of velocity fields and temperatures, respectively. The superscript $k\ge 0$ is used   with $H$ to indicate the Sobolev regularity.} $H=H_1\times H_2=H^0$, $V=V_1\times V_2=H^1$,    and~$U=U_1\times U_2=H^6$.~For any~$T>0$, we~set      
 $$
\XX_T= C(J_T,V)\cap L^2(J_T,H^2)
$$  and  endow this space   with the norm 
$$
\|u\|_{\XX_T}=\|u\|_{C(J_T,V)}+ \|u\|_{L^2(J_T,H^2)}.
$$  
The Leray-type orthogonal   projection onto $H_1$ in $L^2(\mathbb{T}^{3},\R^2)$ is denoted by~$\Pi$.    Applying this   projection   to Eq.~\eqref{0.7}, we eliminate the pressure term and transform   problem~\eqref{0.7},~\eqref{0.8}  into an evolution~system which can be written in the following dimensionless form: 
 \begin{equation}\label{1.3}
 \dot u +Lu+B(u)+Qu=h+\eta,	
 \end{equation}where the unknown is the couple~$u=(v,\te)$, and the linear terms $L$ and $Q$ and the  nonlinear term $B$ are defined by\footnote{With a slight abuse of notation, instead of elements like $B_2(u) $, $Q(u)$,   $Q_1(u)$  with $u=(v,\theta)$, we will often  write  $B_2(v,\theta)$, $Q(v,\theta)$,   $Q_1(v,\theta)$. Note that the term $B_1$ does not depend on $\te$, so we write $B_1(v)$.}
 \begin{align}
   Lu&=(L_1 v, L_2 \te),   \quad   B(u)=(B_1(v), B_2(u) ), \quad  Qu=(Q_1u,0),\nonumber  \\
B_1(v)&=\Pi \left(   \lag v, \nabla\rag v -  \int_0^z \diver v(t,x,y,\z)\,\dd \z\,
\p_z v \right),\nonumber\\ 
B_2(u)&=\lag v,  \nabla\rag \te -  \int_0^z \diver v(t,x,y,\z)\,\dd \z\,
\p_z \te,  \label{E:B2}\\
  Q_1u&= \Pi\left( f   v^\bot  
 -\int_0^z \nabla\te(t,x,y,\z)\,\dd \z\right). \label{E:qban}
\end{align}
  Eq.~\eqref{1.3}  is 
  supplemented with the initial condition
 \begin{equation}\label{1.5}
 u(0)=u_0.	
 \end{equation} 
 \begin{proposition}\label{P:1.1}
For any  $T>0$,   $u_0\in V$,      $\eta\in L^\ty(J_T,H)$, and $h\in H$, there is a   unique solution $u$ of problem~\eqref{1.3},~\eqref{1.5}  belonging  to~$\XX_T$. Let $\RR$ be the mapping  taking the couple\footnote{In what follows, the source term $h$ will be fixed, so we shall  not indicate the dependence   of~$\RR$   on it.} $(u_0,\eta)$   to the solution~$u$.~For any $r>0$, there      is a constant $C=C(r,T)>0$ such that 
$$
\|\RR(u_{0,1}, \eta_1) -\RR( u_{0,2},  \eta_2)\|_{\XX_T} \le 
C \left( \|u_{0,1}- u_{0,2}\|_1+
  \|\eta_1- \eta_2\|_{L^\ty(J_T,H)}\right),  
$$provided that $u_{0,i}\in V$,   $\eta_i \in L^\ty
(J_T,  H)$, and $h\in H$   satisfy   
$$
\|u_{0,i}\|_1+   \| \eta_i\|_{L^\ty(J_T,H)}+\|h\|\le r, \quad i=1,2.
$$\end{proposition}   Existence and uniqueness of  solutions is established in \cite{CT-2007,P-2006, KZ-2007}, and the local Lipschitz property in~\cite{J-2007}.

Inspired by   ideas from   \cite{AS-2005, AS-2006, shirikyan-cmp2006},  together with Eq.~\eqref{1.3}, we  will consider a more general   equation with additional   control     $\zeta$:
  \begin{equation}\label{1.4}
 \dot u +L(u+\zeta)+B(u+\zeta)+Q(u+\zeta)=h+\eta.	
 \end{equation}The   well-posedness of problem \eqref{1.4}, \eqref{1.5}  with $\zeta\in  V$ follows    from that of problem~\eqref{1.3}, \eqref{1.5} using a change of unknown $u'=u+\zeta$. We   denote by~$S(u_0,\zeta,\eta)$ the corresponding solution and by~$S_t(u_0,\zeta,\eta)$ its restriction   at time $t\in J_T$.  To~avoid any ambiguity, in this and next sections, we   write   $S(u_0,0,\eta)$ instead of $S(u_0,\eta)$ defined in Proposition~\ref{P:1.1}. Let  $\pi_1:H\to H_1$ and~$\pi_2:H\to H_2$ be the projections $(v,\te)\mapsto v$ and $(v,\te)\mapsto \te $.
  The following result is proved in Section~\ref{S:6}.

 \begin{proposition}\label{P:1.2} 
 For any  $u_0\in H^4$ and $\zeta, \eta, \xi   \in H^5$ with  $\pi_1 \xi=0 $ and $\pi_2 \zeta =0 $, the following    limits hold in $V$ as $\de\to 0^+$:\footnote{ In \eqref{1.10}, we denote $(0,\Psi)$   the element $\hat \Psi\in H$ such that $\pi_1 \hat \Psi=0$ and $\pi_2 \hat \Psi=\Psi$. Likewise, in what follows, we will often have terms of the form $(v,0)$ that denote an element $\hat v\in H$ with $\pi_1 \hat v=v$ and $\pi_2 \hat v=0$.}
 \begin{gather}
 	\RR_{\de}(u_0,\de^{-\frac12}\zeta,\de^{-1}\eta) \to  u_0+\eta-B(\zeta),\label{1.9}\\
 	 	 	\RR_{\de}(u_0, \delta^{-1} \xi,0) \to  u_0-L\xi-(0,\Psi(u_0,\xi))-Q\xi,\label{1.10}
 \end{gather}          where $\Psi(u_0,\xi)=B_2(\pi_1 u_0-\frac12 Q_1\xi,\pi_2\xi)$.  
   
    \end{proposition}
  
 Now, let $\tilde u=(\tilde v, \tilde \te)= \RR(u_0, 0, \eta)$ be a trajectory of Eq.~\eqref{1.3}   corresponding to    initial condition $u_0\in V$ and   control $\eta \in L^\ty(J_T,H)$.~The  linearization of Eq.~\eqref{1.3} around $\tilde u$ is given by 
  \begin{equation}
 \dot w +Lw+b(\tilde u,w)+Qw=g, \label{1.11}	
  \end{equation}where   $w=(v,\te)$    and  the   term  
  $
  b(\tilde u,w)=(b_1(\tilde v,v),b_2(\tilde u,w)) 
  $  is defined by 
      \begin{align*}
    b_1(\tilde v,v)&=\Pi \bigg(\lag \tilde v, \nabla\rag v+\lag  v, \nabla\rag \tilde v \\&\quad\quad\quad  -  \int_0^z \diver  \tilde v(t,x,y,\z)\,\dd \z\, \p_zv-\int_0^z \diver  v(t,x,y,\z)\,\dd \z\, \p_z \tilde v \bigg), \\
    b_2(\tilde u,w)&=\lag \tilde v,  \nabla\rag \te+\lag v,  \nabla\rag \tilde \te  \\&\quad\quad\quad-  \int_0^z \diver \tilde v(t,x,y,\z)\,\dd \z\,
\p_z\te-  \int_0^z \diver v(t,x,y,\z)\,\dd \z\,
\p_z\tilde\te.  \end{align*}
 Using standard techniques (e.g., see Chapter~III in~\cite{Tem79}), one shows that, for any $w_0\in V$ and $g\in L^2(J_T,H)$, the linear equation~\eqref{1.11}  has a unique solution~$w\in\XX_T$ issued from $w_0$.

\section{Controllability of  the  nonlinear system}\label{S:2}

\subsection{Saturation property   and the result}\label{S:2.1}

In this section, we formulate  a controllability result for    Eq.~\eqref{1.3} that is a   generalization of Theorem~\hyperlink{A}{A}  given   in the Introduction.~We start by introducing some definitions and notation.  
\begin{definition}\label{D:2.1} Let   $\HH$ be a finite-dimensional subspace of $U$. 
 Eq.~\eqref{1.3} is said to be {\it approximately controllable in $H$}
	by   $\HH$-valued controls   if for any $\e>0$, 	any time $T>0$, any initial point~$u_0\in V$, and any target~$u_1\in H$,   there is a control~$\eta
	\in   L^\ty(J_T, \HH) $ such~that
	\begin{equation}\label{2.1}
	\|\RR_T(u_0, 0, \eta) - u_1\| <\e.
	\end{equation} In a similar way, Eq.~\eqref{1.3}  is said to be {\it approximately controllable in $V$} if inequality \eqref{2.1} holds with respect to the~$H^1$-norm~$\|\cdot\|_1$ and the~target~$u_1$ is   arbitrary in~$V$. 	\end{definition}
 	Let  us assume that  $\HH=\HH_1\times \HH_2$, where $\HH_i\subset U_i$, $i=1,2$ are   finite-dimensional subspaces. We denote by $\FF_1(\HH)$ the largest   subspace of $ U_1$ whose  elements    can be approximated,  within any accuracy with respect to the   $H^1$-norm,   by elements of the form    (cf.~\cite{AS-2005, AS-2006, shirikyan-cmp2006})
\begin{equation}\label{2.2}
Q_1(0,\zeta_0)+\zeta_1-\sum_{i=2}^m B_1(\zeta_i),
\end{equation}where    $m\ge2$,   $ \zeta_0\in \HH_2$, and  $\zeta_1, \ldots,\zeta_m \in
\HH_1 $. As $\HH$ is   finite-dimen\-sional, $Q_1$ is linear,  and~$B_1$~is   bilinear, it is easy to see that~$\FF_1(\HH)$ is well~defined and finite-dimensional.

Let $\FF_2(\HH)$ be the   subspace spanned by all the vectors of the form  
\begin{equation}
\xi_0+  	B_2(Q_1(0,\xi_1),\xi_2),
\end{equation}
where  
    $\xi_0,\xi_1, \xi_2\in \HH_2$ are such that $ B_2(Q_1(0,\xi_1),\xi_2)\in U_2$.~We denote by $\FF(\HH)$ the product~$ \FF_1(\HH)\times\FF_2(\HH)$, and define a~non-decreasing sequence~$\{\HH(j)\}$   of     finite-dimensional subspaces   of $U$~by      
\begin{equation}
\HH(0)= \HH,\quad \HH(j)=\FF(\HH(j-1)),    \quad j \geq
1.\label{2.4}
\end{equation}Let us set
\begin{equation}
 \HH(\infty)=  \bigcup_{j=1}^\infty \HH(j).	\label{2.5}
\end{equation}
 \begin{definition}\label{D:2.2}
A   subspace $\HH\subset U$ is~{\it $H$-saturating} (resp.~{\it $V$-satu\-ra\-ting})   if the following two conditions hold:
\begin{itemize}
\item[(a)]	$\HH =\HH_1\times \HH_2$,  where $\HH_i\subset U_i$, $i=1,2$ are   finite-dimensional subspaces;
\item[(b)]	the vector space $\HH(\ty)$   is dense in\footnote{The reason why we take $H_1\times V_2$ and not the   space~$H$ is explained in    Remark~\ref{R:2.7}.} $H_1\times V_2$ (resp.~in~$V$).
\end{itemize}
\end{definition} 
 We are now ready to formulate   the main result of this section.
  \begin{theorem}\label{T:2.3}
 	  If $\HH \subset U$ is an  $H$-saturating  subspace,   then Eq.~\eqref{1.3} is approximately controllable in $H$  by~$\HH$-valued controls. Moreover,  if $\HH$ is    $V$-saturating,  then the equation is approximately controllable in $V$.
 \end{theorem}
 Examples of $H$ and $V$-saturating subspaces are given   in Section~\ref{S:4}. When the  control acts directly only on the temperature component (i.e., $\HH_1=\{0\}$ in~(a) in Definition~\ref{D:2.2}), we   provide an~$H$-saturating subspace---the ten-dimensional space considered in the Introduction.~In particular, Theorem~\hyperlink{A}{A} is obtained as an immediate consequence of Theorems~\ref{T:2.3} and~\ref{T:4.1}.   
  We do not have an example of $V$-saturating subspace acting through the temperature component only. The example we provide is {\it less~degenerate} and combines few modes  from      velocity and temperature components (see Theorem~\ref{T:4.4}). 
Furthermore, in  the case of the 2D and 3D NS systems, there are  necessary and sufficient conditions on the Fourier modes to use in order to have
  approximate controllability (see~\cite{AS-2005, AS-2006, VN-2015}). It would be interesting to obtain similar precise   description in the case of 3D PEs.

\subsection{Proof of Theorem \ref{T:2.3}}

The proof of Theorem \ref{T:2.3} is divided into three steps.~We first show that the temperature and velocity components   can be separately controlled in small time.~Then we derive    simultaneous controllability of both components in   arbitrary fixed~time.

\subsubsection{Controllability of $\theta$-component}

  Let us set    $\HH_i(j)=\pi_i \HH(j)$ for $j\ge0$ and  $i=1,2$. In this subsection, we prove the following proposition.    
 \begin{proposition}\label{P:2.4} Let $\HH_i\subset U_i$, $i=1,2$ be arbitrary   finite-dimensional subspaces
 and   $\HH =\HH_1\times \HH_2$. For any $u_0\in V$ and~$\eta\in \HH_2(\ty)$, there is a family of controls~$\{\zeta_\tau\}_{\tau>0} \subset     L^\ty(J_1, \HH)$ such that 
 \begin{equation}\label{2.6}
 \RR_\tau(u_0, 0, \zeta_\tau)\to u_0+\hat \eta  \quad\text{in $V$ as $\tau\to 0^+$},
 \end{equation}   where $\hat \eta=(0,\eta)\in H$.
 \end{proposition}
   \begin{proof}We first prove the result in the case~$u_0\in U$. It~suffices to show that   for any~$N\ge0$ and~$\eta\in   \HH_2(N)$, there are controls $\{\zeta_\tau\}\subset     L^\ty(J_1, \HH)$ such that  limit~\eqref{2.6} holds.~We argue by induction on~$N\ge0$.

{\it Step~1.~Base case:~$N=0$.}~Let us check that   limit~\eqref{2.6} holds in the case~$N=0$, i.e., for any $\eta\in   \HH_2$. Indeed,   by limit \eqref{1.9} with $\zeta=0$ and $\eta=\hat \eta$, we~have    
	$$ 	
 \RR_{\de}(u_0,0,\de^{-1}\hat \eta)\to  u_0+\hat \eta \quad\text{in $V$ as $\de\to 0^+$}.
$$Taking $\delta=\tau$,
we obtain the required limit with controls $\zeta_\tau=\tau^{-1}\hat\eta$.

	 	{\it Step~2.~Inductive step.}~We  assume that the limit is proved for $N-1$, and take any  $\eta\in   \HH_2(N)$   of the~form
\begin{equation}\label{2.10}
\eta=  \xi_0+   	B_2(Q_1(0,\xi_1),\xi_2)  
\end{equation}with some       $\xi_0,\xi_1, \xi_2\in \HH_2(N-1)$.    Let  us set $\hat\xi_i= (0,\xi_i)\in H$, $i=1,2,3$. 
   Using limit~\eqref{1.10} with~$\xi=\hat \xi_1$, we get
\begin{equation}\label{2.11}
 \RR_{\de}(u_0, \delta^{-1} \hat\xi_1,0) \to  u_0-L\hat\xi_1-(0,\Psi(u_0,\hat\xi_1))-Q\hat\xi_1.
\end{equation}
By the uniqueness of the solution of the Cauchy problem,
 the following  equality holds for any $t\ge0$:
\begin{equation}\label{2.12}
\RR_{t}(u_0+\delta^{-1} \hat\xi_1,0,0)=\RR_{t}(u_0,\delta^{-1} \hat\xi_1,0)+\delta^{-1} \hat\xi_1. 
\end{equation}Taking here $t=\delta$   and using \eqref{2.11}, we obtain the limit
\begin{equation}\label{2.12b}
\|\RR_{\delta}(u_0+\delta^{-\frac12}\hat \xi_1,0,0)-u_0 +L\hat\xi_1+(0,\Psi(u_0,\hat\xi_1))+Q\hat\xi_1    -\delta^{-1} \hat\xi_1\|_1\to 0		
\end{equation}as $\delta\to 0^+.$   The fact that $\xi_1 \in   \HH_2(N-1)$ and  the   induction   hypothesis  imply that, for any $\de>0$, there is a family of controls       $\{\zeta_{\tau,\de}^1\} \subset   L^\ty(J_1, \HH)$ such that
$$
	 \RR_\tau(u_0, 0, \zeta_{\tau,\de}^1)\to  u_0+\delta^{-\frac12}\hat \xi_1   \quad\text{in $V$ as $\tau\to 0^+$}.	
$$
From Proposition~\ref{P:1.1} it follows   
$$
	 \RR_\de(\RR_\tau(u_0, 0, \zeta_{\tau,\de}^1),0,0)\to   \RR_\de(u_0+\delta^{-\frac12}\hat \xi_1,0,0)   \quad\text{in $V$ as $\tau\to 0^+$}.	
$$Combining this with \eqref{2.12b}, we find a family   $\{\zeta_\de^2\} \subset   L^\ty(J_1, \HH)$ verifying 
$$
\|\RR_{\delta}(u_0,0,\zeta_\de^2)-u_0 +L\hat\xi_1+(0,\Psi(u_0,\hat\xi_1))+Q\hat\xi_1    -\delta^{-1} \hat\xi_1\|_1\to 0		
$$as $\delta\to 0^+.$ Using one more time     the assumption   $\xi_1 \in   \HH_2(N-1)$,  the   induction   hypothesis, and Proposition~\ref{P:1.1}, we find       $\{\zeta_\tau^3\} \subset   L^\ty(J_1, \HH)$ such~that 
\begin{equation}\label{2.13}
	 \RR_\tau(u_0, 0, \zeta_\tau^3)\to u_0- L\hat\xi_1-(0,\Psi(u_0,\hat\xi_1)) -Q\hat\xi_1  \quad\text{in $V$ as $\tau\to 0^+$}.	
\end{equation}
  Now we use the following lemma.
\begin{lemma}\label{L:2.5b}
	Let us denote 
$$
F_{\xi_1}(u_0)=u_0- L\hat\xi_1-(0,\Psi(u_0,\hat\xi_1)) -Q\hat\xi_1. 
$$Then  	
\begin{equation}\label{E:2.12b}
  F_{-\xi_2}(F_{-\xi_1}(F_{\xi_2}(F_{\xi_1}(u_0))))=u_0+(0, B_2(Q_1(0,\xi_1-\xi_2),\xi_1+\xi_2)).
\end{equation}
\end{lemma}
Using \eqref{E:2.12b} and
 iterating four times the argument of the   construction of   the family    $\{\zeta_\tau^3\}$, we   find a family $\{\zeta_\tau^4\} \subset   L^\ty(J_1, \HH)$ such that 
$$
	 \RR_\tau(u_0, 0, \zeta_\tau^4)\to u_0+(0, {B_2(Q_1(0,\xi_1-\xi_2),\xi_1+\xi_2)})  \quad\text{in $V$ as $\tau\to 0^+$}.	
$$ 
  This and the induction   hypothesis imply that, for any   $\xi_0,\xi_1, \xi_2\in \HH_2(N-1)$,  there are  
controls     $\{\zeta_\tau^5\} \subset   L^\ty(J_1, \HH)$ such~that  
$$
	 \RR_\tau(u_0, 0, \zeta_\tau^5)\to u_0+(0,\xi_0+ B_2(Q_1(0,\xi_1),\xi_2))    \quad\text{in $V$ as $\tau\to 0^+$}.	
$$Iterating this argument, we show that    the system can be controlled in small time  to any target of the form $u_0+\hat \eta$ (in the sense of limit \eqref{2.6}), where $\hat\eta$ is now a linear combination of vectors of the form \eqref{2.10}. 
 This completes the proof of the proposition in the case of a regular initial condition $u_0$.
 In the case~$u_0\in V$,  it suffices to take control equal to zero on a   small time~interval, to use     the  regularizing property of the PEs (e.g., see Theorem~3.1 in~\cite{P-2006}), and   apply the already proved result for  regular  initial condition.
 \end{proof}
 \begin{proof}[Proof of Lemma \ref{L:2.5b}] For any   smooth $u$ in $H$, we have
\begin{align}\label{E:dzb}
F_{-\xi_1}(F_{\xi_2}(u))&= F_{\xi_2}(u)+ L\hat\xi_1+(0,\Psi(F_{\xi_2}(u),\hat\xi_1)) +Q\hat\xi_1 \nonumber \\
&=u+ L(\hat\xi_1-\hat\xi_2)-(0,\Psi(u,\hat\xi_2))+(0,\Psi(F_{\xi_2}(u),\hat\xi_1)) +Q(\hat\xi_1-\hat\xi_2).
\end{align}
Replacing in this equality $u$ by $F_{\xi_1}(u_0)$, we obtain
\begin{align}\label{54597}
F_{-\xi_1}(F_{\xi_2}&(F_{\xi_1}(u_0))) =u_0- L \hat\xi_2-Q \hat\xi_2\nonumber\\&\quad-(0,\Psi(u_0,\hat\xi_1))-(0,\Psi(F_{\xi_1}(u_0),\hat\xi_2))+(0,\Psi(F_{\xi_2}(F_{\xi_1}(u_0)),\hat\xi_1)) .
\end{align}Note that 
\begin{equation}\label{5564}
\Psi(u,\xi)  \text{ does not depend on $\pi_2u$ and $\pi_1 L\hat \xi_1=\pi_1 L\hat \xi_2=0$,}
\end{equation}  so using \eqref{E:dzb}, we get
\begin{align*}
\Psi(F_{\xi_1}(u_0),\hat\xi_2)&=\Psi (u_0-Q\hat\xi_1 ,\hat\xi_2), \\	
 \Psi(F_{\xi_2}(F_{\xi_1}(u_0)),\hat\xi_1)&=  \Psi( u_0+ Q(\hat\xi_1-\hat\xi_2),\hat\xi_1).
\end{align*}Thus \eqref{54597} can be rewritten as
\begin{align}
F_{-\xi_1}(F_{\xi_2}&(F_{\xi_1}(u_0))) =u_0- L \hat\xi_2-Q \hat\xi_2\nonumber\\&\quad+(0, \Psi( u_0+ Q(\hat\xi_1-\hat\xi_2),\hat\xi_1)-\Psi (u_0-Q\hat\xi_1 ,\hat\xi_2)-\Psi(u_0,\hat\xi_1)), \nonumber
\end{align}hence
\begin{align*}
 F_{-\xi_2}(F_{-\xi_1}&(F_{\xi_2}(F_{\xi_1}(u_0)))) =u_0 + (0,\Psi ( F_{-\xi_1}(F_{\xi_2} (F_{\xi_1}(u_0))) , \hat \xi_2) ) \nonumber\\&\quad+(0, \Psi( u_0+ Q(\hat\xi_1-\hat\xi_2),\hat\xi_1)-\Psi (u_0-Q\hat\xi_1 ,\hat\xi_2)-\Psi(u_0,\hat\xi_1))\nonumber\\
 & =u_0 + (0,\Psi (  u_0-Q\hat\xi_2,\hat\xi_2)) \nonumber\\&\quad+(0, \Psi( u_0+ Q(\hat\xi_1-\hat\xi_2),\hat\xi_1)-\Psi (u_0-Q\hat\xi_1 ,\hat\xi_2)-\Psi(u_0,\hat\xi_1))\nonumber\\&=u_0+(0,\mathfrak{b}_2(\xi_1-\xi_2,\xi_1+\xi_2)),
\end{align*}where we used again \eqref{5564} and the equality
$$
\Psi ( F_{-\xi_1}(F_{\xi_2} (F_{\xi_1}(u_0))) , \hat \xi_2)=\Psi (u_0-Q \hat\xi_2,\hat\xi_2).
$$ \end{proof}

\subsubsection{Controllability of $v$-component}

Here we prove the following  version of Proposition~\ref{P:2.4} for the $v$-component.  
 \begin{proposition}\label{P:2.7}
 Let $\HH_i\subset U_i$, $i=1,2$ be     finite-dimensional subspaces,  let $\HH =\HH_1\times \HH_2$, and assume that $\HH_2(\ty)$ is dense in $V_2$.   For any~$u_0\in V$ and~$\eta\in \HH_1(\ty)$,  there is a family of    controls     $\{\xi_\tau\}_{\tau>0} \subset    L^\ty(J_1, \HH)$ such that  
  \begin{equation}\label{2:11}
 \RR_\tau(u_0, 0, \xi_\tau)\to u_0+\hat \eta  \quad\text{in $V$ as $\tau\to 0^+$},
 \end{equation}   where $\hat \eta=(\eta,0)\in H$.
 \end{proposition}
 \begin{proof}
 The argument is   close to the one used in Proposition~\ref{P:2.4}. Again, without loss of generality, we can assume that $u_0\in U$. We prove limit \eqref{2:11} for any~$\eta\in \HH_1(N)$,  arguing by induction on~$N\ge0$. The base case~$N=0$ follows from limit~\eqref{1.9} with~$\zeta=0$ and~$\eta= \hat \eta$:
 $$ 	
 \RR_{\de}(u_0,0,\de^{-1}\hat\eta)\to  u_0+\hat\eta \quad\text{in $V$ as $\de\to 0^+$}.
$$Taking $\delta=\tau$,
we obtain the required limit with $\xi_\tau=\tau^{-1}\hat\eta$.

Assume that the limit is proved in the case $N-1$, and let $\eta\in \HH_1(N)$.    
By   approximation, we can suppose  that $\eta$ is of the form
$$
\eta= Q_1\hat \zeta_0 +\zeta_1  -\sum_{i=2}^m B_1(\zeta_i) 
$$for some       $m\ge2$,    $\hat \zeta_0=(0,\zeta_0)$, $\zeta_0 \in \HH_2(N-1)$, and $ \zeta_1, \ldots,\zeta_m \in\HH_1(N-1) $.

{\it Step~1.~Direction $\zeta_1  -\sum_{i=2}^m B_1(\zeta_i)$.}
   Limit \eqref{1.9} with   $\zeta=\hat \zeta_2=(\zeta_2,0)$ and~$\eta=0$ implies that  \begin{equation}\label{2:13}
 \RR_{\de}(u_0,\de^{-\frac12}\hat\zeta_2,0)\to  u_0- B(\hat \zeta_2)  \quad\text{in $V$ as $\de\to 0^+$}.
\end{equation}
 The equality
$$
\RR_{\delta}(u_0+\delta^{-\frac12}\hat \zeta_2,0,0)=\RR_{\delta}(u_0,\delta^{-\frac12}\hat \zeta_2,0)+\delta^{-\frac12}\hat \zeta_2 
$$ and  limit \eqref{2:13}  show  that
$$
\|\RR_{\delta}(u_0+\delta^{-\frac12}\hat \zeta_2,0,0)-u_0+B(\hat \zeta_2)-\delta^{-\frac12}\hat \zeta_2\|_1\to 0\quad \text{as $\delta\to 0^+.$}	
$$   
  The fact that $   \zeta_2\in  \HH_1(N-1)$   and  the   induction   hypothesis  imply that, for any $\de>0$, there is a family of controls       $\{\xi_{\tau,\de}^1\} \subset   L^\ty(J_1, \HH)$ such that
$$
	 \RR_\tau(u_0, 0, \xi_{\tau,\de}^1)\to  u_0+\delta^{-\frac12}\hat \zeta_2   \quad\text{in $V$ as $\tau\to 0^+$}.	
$$
Then Proposition~\ref{P:1.1} implies that   
$$
	 \RR_\de(\RR_\tau(u_0, 0, \xi_{\tau,\de}^1),0,0)\to   \RR_\de(u_0+\delta^{-\frac12}\hat \zeta_2,0,0)   \quad\text{in $V$ as $\tau\to 0^+$}.	
$$Combining this with \eqref{2:13}, we find a family   $\{\xi_\de^2\} \subset   L^\ty(J_1, \HH)$ verifying 
$$
\|\RR_{\delta}(u_0,0,\xi_\de^2)-u_0+B(\hat \zeta_2)-\delta^{-\frac12}\hat \zeta_2\|_1\to 0		
$$as $\delta\to 0^+.$ Using        the assumption     $ \zeta_1, \zeta_2\in  \HH_1(N-1)$, the   induction   hypothesis, and Proposition~\ref{P:1.1},   we find  a family of        controls     $\{\xi_\tau^3\}\subset    L^\ty(J_1, \HH)$ such~that  
$$
 \|\RR_\tau(u_0, 0, \xi_\tau^3)- u_0-(\zeta_1-B_1(\zeta_2), 0)\|_1\to 0  \quad\text{ as $\tau\to 0^+$}.	
$$
 Iterating this argument   with     $\zeta_3, \ldots, \zeta_m$, we construct a family $\{\xi_\tau^{m+1}\}\subset    L^\ty(J_1, \HH)$ such that
 \begin{equation}\label{2:14}
 \|\RR_\tau(u_0, 0, \xi_\tau^{m+1})-  u_0-(\zeta_1-\sum_{i=2}^m B_1(\zeta_i),0)\|_1 \to 0       \quad\text{ as $\tau\to 0^+$}.	
 \end{equation}

{\it Step~2.~Direction $Q_1\hat\zeta_0$.} 
Let~$\hat u_0\in H^4$.
 By limit~\eqref{1.10} with $\xi=\hat \zeta_0$ and $\eta=0$, we have 
\begin{gather*}
 	 	 	\RR_{\de}(\hat u_0, \delta^{-1} \hat\zeta_0,0) \to  \hat u_0-L\hat\zeta_0-(0,\Psi( \hat u_0, {\hat\zeta}_0))-Q\hat\zeta_0\\=\hat u_0 -(Q_1\hat\zeta_0, L_2\zeta_0+\Psi( \hat u_0, {\hat\zeta}_0)). 
\end{gather*}
The equality
$$
\RR_{\de}(\hat u_0+\delta^{-1} \hat\zeta_0,0 ,0)=\RR_{\de}(\hat u_0, \delta^{-1} \hat\zeta_0,0)+\delta^{-1} \hat\zeta_0
$$implies that 
$$
\|\RR_{\de}(\hat u_0+\delta^{-1} \hat\zeta_0,0 ,0) - \hat u_0+(Q_1\hat\zeta_0, L_2\zeta_0+\Psi(\hat u_0, {\hat\zeta}_0))- \delta^{-1} \hat\zeta_0 \|_1\to 0
$$
as $\delta\to0^+$.~Combining this   with the assumption   that~$\HH_2(\ty)$ is dense in $V_2$ and Propositions~\ref{P:1.1} and~\ref{P:2.4}, we construct a family of controls   $ \{\xi_\tau^{m+2}\}\subset     L^\ty(J_1, \HH)$ such~that  
$$
 \RR_\tau(\hat u_0,0,  \xi_\tau^{m+2})\to \hat  u_0-(Q_1\hat\zeta_0, 0) 
 \quad\text{in $V$ as $\tau\to 0^+$}.	
$$ 
Taking   
$$
\hat u_0=  u_0+(\zeta_1-\sum_{i=2}^m B_1(\zeta_i), 0)   
$$
and using \eqref{2:14}, we find a family of    controls
  $ \{\xi_\tau \}\subset     L^\ty(J_1, \HH) $ such that limit~\eqref{2:11}~holds.
  \end{proof}

\subsubsection{Completion of the proof}\label{S:2.2.3}

   Assume that $\HH\subset U$ is an  $H$-saturating (resp.~$V$-saturating) subspace, and  let~$\e>0$, $T>0$, $u_0\in V$, and   $u_1\in H$ (resp.~$u_1\in V$) be arbitrary. Then
   there is~$\eta=(\eta_1,\eta_2) \in \HH(\ty)$   such that 
\begin{equation}\label{EEq}
    \|\RR_T(u_0, 0, 0)+\eta-   u_1\| <\frac\e2 \quad \left(\text{resp.~}  \|\RR_T(u_0, 0, 0)+\eta-  u_1\|_1 < \frac\e2 \right).   	
\end{equation}  Let us denote $\hat u_0=\RR_T(u_0, 0, 0)+\eta$ and take   $t_0>0$ and $r>0$ so small that 
\begin{equation}\label{2:15}
    \|\RR_{t}(u, 0, 0) - \hat u_0 \|_1 <\frac\e2 \quad \text{for $t\in [0,t_0]$ and $u\in B_V(\hat u_0,r)$}.   
 \end{equation}This is possible by   Proposition~\ref{P:1.1}. Choosing, if necessary, $t_0$ smaller, we will also have 
   \begin{equation} \label{2:16}
    \|\RR_{T-t}(u_0, 0, 0) -\RR_{T}(u_0, 0, 0)\|_1<\frac r2	\quad \text{for $t\in [0,t_0]$}.
   \end{equation} Now applying Propositions~\ref{P:2.4} and~\ref{P:2.7}   with initial condition $S_{T-t_0}(u_0, 0, 0)$,   we find  a time $\tau \in (0,t_0)$ and a control  $\xi\in L^\ty([0,\tau], \HH)$ such that 
$$
   \|S_\tau(S_{T-t_0}(u_0, 0, 0), 0,  \xi)-S_{T-t_0}(u_0, 0, 0)-(\eta_1, \eta_2)\|_1<\frac r2.
$$In view of \eqref{2:16}, this implies 
  that $ S_\tau(S_{T-t_0}(u_0, 0, 0),  0, \xi)\in B_V(\hat u_0, r). $ 
 Finally, using~\eqref{EEq} and \eqref{2:15}, we   conclude that 
   $$
  \|\RR_T(u_0, 0, \zeta)  -   u_1\| <\e \quad \left(\text{resp.~}  \|\RR_T(u_0, 0,\zeta) -  u_1\|_1 < \e \right),  	
$$where 
 $
\zeta(t)= \I_{[T-t_0,T-t_0+\tau]}  \xi(t-T+t_0)$, $t\in J_{T}. 
$   This completes the proof of~Theorem~\ref{T:2.3}.
\begin{remark} 
Note that the above proof gives approximate controllability to any target~$u_1$ in~$H_1\times V_2$ with respect to the norm of that space.
 \end{remark}
\begin{remark}\label{R:2.7}
The assumption that $\HH_2(\ty)$ is dense in $V_2$ (see (b) in Definition~\ref{D:2.2}) plays an important role in the above proof. We use  it in Step~2 of the proof of Proposition~\ref{P:2.7}. If $\HH_2(\ty)$   was  dense only in $H$, we would need a version of Proposition~\ref{P:1.1}   with respect to the 
$L^2$-norm. The latter is an open problem.
\end{remark}

 \section{Controllability of   linearized system}\label{S:3}

  \subsection{Saturation   for  linearized system  and the result}\label{S:3.1}

Before formulating the main result of this section, let us define a saturation property for    linearized system \eqref{1.11},  which is different from the one   used in the        nonlinear case (cf. Definition~\ref{D:2.2}),    and recall  the concept  of    {\it observable measures} introduced in~\cite{KNS-2018}.

 We   assume that  $\HH=\HH_1\times \HH_2$, where $\HH_1=\{0\}\subset U_1$ and 
 $\HH_2\subset U_2$  is  a finite-dimensional subspace. Let us   define 
 vector spaces $\GG_1(\ty)\subset U_1$ and~$\GG_2(\ty)\subset U_2$  as follows:
\begin{description}
	\item [$\bullet$] $\GG_2(\ty)=\cup_{j=0}^\ty\GG_2(j)$, where
	  $\GG_2(0)=\HH_2$ and $\GG_2(j),$~$j\ge1$ is      the   space spanned by all the vectors of the~form  
$$
\xi_0+ \mathfrak{b}_2(\xi_1,\xi_2),
$$where   $\xi_0, \xi_1\in \GG_2(j-1)$ and $\xi_2\in \HH_2$   are such that 
$$
\mathfrak{b}_2(\xi_1,\xi_2)=B_2(Q_1(0,\xi_1),\xi_2)-B_2(Q_1(0,\xi_2),\xi_1)\in U_2;
$$
 \item [$\bullet$]  $\GG_1(\ty)$   is      the   space spanned by all the vectors of the~form 
$$
 Q_1(0,\zeta_0)+ b_1(Q_1(0,\zeta_1),Q_1(0,\zeta_2)),
$$where  
 $\zeta_0, \zeta_1\in \GG_2(\ty)$ and $\zeta_2\in  \HH_2$ are such that  $$Q_1(0,\zeta_0)+ b_1(Q_1(0,\zeta_1),Q_1(0,\zeta_2))\in U_1.$$
\end{description}We set $\GG(\ty)=\GG_1(\ty)\times \GG_2(\ty)$.

  \begin{definition}\label{D:3.1}
  A   finite-dimensional  space $\HH$ as above is said to be   {\it  saturating}  for     linearized system~\eqref{1.11} if~$\GG(\infty)$ is dense in $H$.
\end{definition} 
We will see  in Section~\ref{S:4} that   the ten-dimensional subspace  defined in the Introduction is saturating in the sense of Definition~\ref{D:3.1}.  
Let us take any~$T>0$, denote
  $\EE=L^\ty(J_T,\HH)$, and let     $\{\varphi_i\}_{i=1}^d$ be a basis in $\HH$. 
 \begin{definition} \label{D:3.2}
   A function~$\zeta\in \EE$ is said to be   {\it observable}    if for any continuously differentiable functions $a_i:J_T\to\R$, $i\in[\![1,d]\!]$ and any continuous function $a_0:J_T\to\R$  	the equality 
$$
		\sum_{i=1}^d a_i(t)\lag\zeta(t),\varphi_i\rag-a_0(t)=0\quad\mbox{ for almost every $t\in J_T$}
$$
	implies that $a_i(t)=0$ for any $t\in J_T$ and  $i\in[\![0,d]\!]$. 
	 A   measure~$\ell\in \PP(\EE)$   is said to be   observable    if~$\ell$-almost every trajectory $\eta\in \EE$ is observable. 
\end{definition}
It is easy to see that  the    observability does not depend on the choice of the basis~$\{\varphi_i\} $ in $\HH$. See Section~5 in~\cite{KNS-2018} for examples of 
observable measures. In particular, it is shown there that the 
  law of the Haar   noise defined by~\eqref{0.11},~\eqref{0.12} is    observable.

Let $u\in U$, and  let~$D_\eta \RR_T(u,\eta)$ be the  derivative of $\RR_T(u,\eta)$ with respect to~$\eta\in \EE$.~Then the linear~mapping
$$
D_\eta \RR_T(u,\eta): \EE\to U, \quad  g\mapsto w(T)
$$ is the resolving operator for Eq.~\eqref{1.11}, where $\tilde u (t)= S_T(u,\eta)$, $(u,\eta)\in U\times \EE$, and~$t\in J_T$.  Let $\KK^u$ be a Borel set in~$\EE$ defined by 
\begin{equation}\label{E:ku}
\KK^u=\{\eta\in \EE:  \text{the image of~$D_\eta \RR_T(u,\eta)$ is dense in $U$}\}.
\end{equation}  
    
   The following  theorem   can be seen as a non-Gaussian extension of the non-degeneracy
     property of the Malliavin matrix. The latter
     is known to be an important ingredient in the  study of ergodicity and existence of positive densities for stochastic equations driven by a white-in-time noise (see~\cite{MP-2006, HM-2006, HM-2011, FGRT-2015, HM-2015}).

 \begin{theorem}\label{T:3.3}
	Let~$\ell \in \PP(\EE)$, and let $\HH $ be a  saturating subspace   in the sense of Definition~\ref{D:3.1}.    If there is $\tau\in (0,T)$ such that the restriction\footnote{This means that
	$\ell'$  is the image of $\ell$ by the mapping   $\pi_{J_T}:\EE\to L^\ty(J_\tau,\HH),$   $\eta\mapsto \eta|_{J_\tau}$.} $\ell'$ of~$\ell$ to the interval~$J_\tau$    is observable, then~$\ell(\KK^u)=1$ for any $  u\in U$.
	\end{theorem}
 In other words, the conclusion of this theorem is  that  Eq.~\eqref{1.11} is approximately controllable in $V$ by $\HH$-valued control $g$ for any~$ u\in U$ and~$\ell$-a.e.~$\eta\in \EE$. 
 
\subsection{Proof of Theorem~\ref{T:3.3}}

We follow the    scheme    used in the case of the      complex Ginzburg--Landau equation~considered in \cite{KNS-2018}.~Let~$w(t;w_0,g)$ be the   solution of~ Eq.~\eqref{1.11} corresponding~to initial condition~$w_0 \in H$, control $g\in  \EE$, and reference trajectory $\tilde u (t)= S_t(u,\eta)$, $t\in J_T$.~Our goal is to prove that the vector space 
$\La=\{w(T;0,g),g\in \EE\}$  is dense in~$U$
   for any~$ u\in U
    $ and~$\ell$-a.e.~$\eta\in \EE$.

   A well-known property of approximate controllability
  by   initial condition,\footnote{In the case of Eq.~\eqref{1.11}, this   can be proved by literally repeating the arguments of   Section~7.2 in~\cite{KNS-2018}, where a similar result is proved for linear parabolic equations.} applied to Eq.~\eqref{1.11} with  $g\equiv0$, shows that the vector space
  $\{w(s;w_0,0),w_0\in H\}$
    is dense in~$U$ for any $s\in [0,T]$. 
  Let us   apply this   result for the interval  $[\tau,T]$, where~$\tau$ is as in Theorem~\ref{T:3.3}.  
   Furthermore,   the resolving operator for Eq.~\eqref{1.11} on the interval $[\tau, T]$ with $g\equiv0$  is continuous  from $H$ to $U$. Hence, to show that $ \La$ is dense in~$U$, it suffices to prove the density of the vector space~$\{w(\tau;0,g),g\in \EE\}$  in $H$. 
  	
For any  $0\le s\le t\le \tau$, we denote by    $R^{\tilde u}(t,s):H\to H$   the two-parameter resolving operator for the homogeneous problem
\begin{equation}\label{3.3}
	 \dot w +Lw+b(\tilde u,w)+Qw=0,	\quad w(s)=w_0.
\end{equation}
Let $G^{\tilde u}$ be the
{\it controllability Gramian} for  Eq.~\eqref{1.11} (see~Chapter~1 in~\cite{coron2007}):
$$
	G^{\tilde u}=\int_0^\tau R^{\tilde u}(\tau,t)\ppP_\HH R^{\tilde u}(\tau,t)^*\dd t,
$$
	where $R^{\tilde u}(\tau,t)^*:H\to H$ is the adjoint of $R^{\tilde u}(\tau,t)$, and~$\ppP_\HH $ is   the orthogonal projection onto~$\HH$ in $H$.  It is easy to see that the required assertion will
be established if we show that $\Ker G^{\tilde u}$ is trivial for $\ell$-a.e. $\eta\in \EE$.

 	It is easily seen that 		$p(t)=R^{\tilde u}(\tau,t)^*w_0$
  is the solution of the dual of problem~\eqref{3.3} given by
	  \begin{equation}
 \dot p -L p-b(\tilde u)^*p-Q^*p=0, \quad  p(\tau)=w_0, \label{3.4}	
  \end{equation}where $b(\tilde u)^*$ and $Q^*$ are the   adjoints %\footnote{We do not give the formulas of $b(\tilde u)^*$ and $Q^*$, since we will not use them.}  
  of $b(\tilde u, \cdot)$ and $Q$ in $H$.

Let us fix any observable $\eta\in \EE$  and show that $\Ker(G^{\tilde u})=\{0\}$. For any    $w_0\in \Ker(G^{\tilde u})$, we have 
$$
\lag G^{\tilde u}w_0,w_0\rag=\int_0^\tau\|\ppP_\HH R^{\tilde u}(\tau,t)^*w_0\|^2\dd t=\int_0^\tau\|\ppP_\HH p(t)\|^2\dd t=0,
$$which implies that 
  $\ppP_\HH p(t)=0$ for      $t\in J_\tau$. Thus,  
\begin{equation} \label{3.5}
	\lag \zeta,p(t)\rag=0,\quad \text{  $t\in J_\tau$}
\end{equation}for any   $\zeta\in\HH$.  From this we are going to derive   that $\pi_1 w_0$ and $\pi_2 w_0$ are zero.
   
 {\it Step~1. Proof of $\pi_2 w_0=0$.} Let us denote $p_i(t)=\pi_i p(t)$, $i=1,2$. In this step,  we show that
 \begin{equation}\label{EE:1a}
 	     p_2(t)=0\quad  \text{for  $t\in J_\tau$}.
 \end{equation}Choosing $t=\tau$ in this equality, we    get  $\pi_2 w_0=0$. To prove~\eqref{EE:1a}, let us take   $\zeta=\hat\xi=(0,\xi) \in \{0\}\times \HH_2$ in   \eqref{3.5}. Then\footnote{We shall use the same notation  $\lag \cdot, \cdot \rag$  for   the scalar products in the spaces $H$, $H_1$, and $H_2$.} 
 \begin{equation} \label{EE:2a}
	\lag \xi,p_2(t)\rag=0,\quad \text{  $t\in J_\tau$}.
\end{equation}  This shows that $p_2(t)$ is orthogonal to $\HH_2$ for any $t\in J_\tau$. 
 In what follows, we prove that~$p_2(t)$ is orthogonal to all subspaces~$\GG_2(j), j\ge1$. By the saturation assumption, the subspace    $\GG_2(\ty)$ is dense in $H_2$, so we   get~\eqref{EE:1a}.

  We proceed by induction on~$j\ge0$. 
  The case $j=0$ is already considered above.
Assuming  that   \eqref{EE:2a} holds for~any $\xi\in\GG_2(j-1)$, let us prove it for any~$\xi\in\GG_2(j)$. Taking  $\zeta=\hat\xi=(0,\xi) \in \{0\}\times \HH_2$ in~\eqref{3.4},  differentiating the resulting equality  in time, and using~\eqref{3.4}, we obtain  
\begin{equation}\label{EE:RET}	
\lag  Q_1 \hat \xi,p_1(t)\rag + \lag L_2\xi+b_2(\tilde u(t),\hat \xi), p_2(t)\rag=0,\quad\text{ $t\in J_\tau$}.
\end{equation}
Note that, as $\pi_1\hat \xi=0$, we have $b_2(\tilde u(t),\hat \xi)=B_2(\tilde v(t), \xi)$, where  $\tilde v(t)=\pi_1 \tilde u(t)$. Thus \eqref{EE:RET} becomes
$$\lag  Q_1 \hat \xi,p_1(t)\rag +\left \lag L_2\xi+B_2(\tilde v(t), \xi), p_2(t)\right\rag=0,\quad\text{ $t\in J_\tau$}.
$$
  Taking the derivative  in time of this equality, we get
$$	\lag  Q_1 \hat \xi,\dot p_1(t)\rag+
\left \lag  B_2(\dot {\tilde v}(t),\xi), p_2(t)\right\rag+\left \lag L_2\xi+B_2(\tilde v(t),\xi), \dot p_2(t)\right\rag=0,\quad\text{ $t\in J_\tau$}.
$$
From the equations for $\tilde v$ and $p$ and the fact that $\pi_1\eta=0$    we derive  
\begin{gather}
  \lag L_1 (Q_1 \hat \xi) +b_1(\tilde v(t), Q_1 \hat \xi) + Q_1 \hat q(t) ,p_1(t)\rag +\left \lag L_2q(t)+b_2(\tilde u(t), \hat q(t)),  p_2(t)\right\rag
\nonumber\\
	- \left \lag  B_2(  L_1\tilde v(t)+B_1(\tilde v(t))+Q_1  \tilde u(t)-h_1 ,\xi), p_2(t)\right\rag =0,  \label{EE:3a}
\end{gather}where 
$q(t)= L_2\xi+B_2(\tilde v(t),\xi)$ and $\hat q(t)=(Q_1 \hat \xi,q(t))\in H$,  $t\in J_\tau$.  
   Setting
\begin{equation}\label{3.7}
	y(t)=\tilde u(t)-\int_0^t\eta(s)\,\dd s=\tilde u(t)-\sum_{i=1}^d \varphi_i \int_0^t\eta^i(s)\,\dd s,	
\end{equation}
where $\eta^i(t)=\lag\eta(t),\varphi_i\rag$, and using the equalities 
\begin{align}\label{E:ffgdh}
b_2(\tilde u(t), \hat q(t))&=B_2(\tilde v(t),  q(t))+B_2( Q_1\hat\xi, \pi_2 \tilde u(t))\nonumber\\&=B_2(\tilde v(t),  q(t))+B_2(Q_1\hat\xi, \pi_2 y(t) )+\sum_{i=1}^d B_2(Q_1\hat\xi, \pi_2\varphi_i)\int_0^t\eta^i(s)\,\dd s,
\end{align}
 we rewrite~\eqref{EE:3a}   as 
\begin{gather*}
  \lag L_1 (Q_1 \hat \xi) +b_1(\tilde v(t), Q_1 \hat \xi) + Q_1 \hat q(t) ,p_1(t)\rag \\+ \lag L_2q(t)+B_2(\tilde v(t),  q(t))+B_2(Q_1\hat\xi, \pi_2 y(t) ),  p_2(t)\rag
\nonumber\\
	-  \lag  B_2(  L_1\tilde v(t)+B_1(\tilde v(t))+Q_1  y(t)-h_1 ,\xi), p_2(t)\rag\nonumber\\+
	\sum_{i=1}^d \lag  B_2(Q_1\hat\xi,\pi_2\varphi_i)- B_2(Q_1   \varphi_i ,\xi), p_2(t)\rag
	  \int_0^t\eta^i(s)\,\dd s =0.  
\end{gather*}
Taking the derivative   in time of this equality and setting 
\begin{align*}
	a_i(t)&=  \lag  B_2(Q_1\hat\xi,\pi_2\varphi_i)- B_2(Q_1   \varphi_i ,\xi), p_2(t)\rag, \quad i\in[\![1,d]\!], \\
	a_0(t)&=\frac{\dd}{\dd t} \Big(   \lag L_1 (Q_1 \hat \xi) +b_1(\tilde v(t), Q_1 \hat \xi) + Q_1 \hat q(t) ,p_1(t)\rag\\&\quad\quad + \lag L_2q(t)+B_2(\tilde v(t),  q(t))+B_2(Q_1\hat\xi, \pi_2 y(t) ),  p_2(t)\rag
 \\&\quad\quad -  \lag  B_2(  L_1\tilde v(t)+B_1(\tilde v(t))+Q_1  y(t)-h_1 ,\xi), p_2(t)\rag\Big)\\&\quad +
	\sum_{i=1}^d \lag  B_2(Q_1\hat\xi,\pi_2\varphi_i)- B_2(Q_1   \varphi_i ,\xi), \dot p_2(t)\rag
	  \int_0^t\eta^i(s)\,\dd s , 
\end{align*}
we obtain 
$$
a_0(t)+\sum_{i=1}^d a_i(t)\eta^i(t)=0. 
$$
The functions~$\{a_i\}_{i=1}^d $ are continuously differentiable and $a_0$ is continuous.~The observability of~$\eta$ implies  that $a_i(t)= 0$ for $t\in J_\tau$ and  $i\in[\![0,d]\!]$. Thus \eqref{EE:2a} holds with $\xi $ replaced by   
\begin{equation}\label{ERTE}
 B_2(Q_1\hat\xi,\pi_2\varphi_i)- B_2(Q_1   \varphi_i ,\xi)=\mathfrak{b}_2(\xi,\pi_2\varphi_i).
\end{equation}We conclude that \eqref{EE:2a} holds with any $\xi    \in \GG_2(j)$.

  {\it Step~2. Proof of $\pi_1 w_0=0$.} 
  In this step,  we show that
 \begin{equation}\label{EE:1b}
 	     p_1(t)=0\quad  \text{for  $t\in J_\tau$}.
 \end{equation}Choosing $t=\tau$, we    get  $\pi_1 w_0=0$, which will complete the proof of the theorem.

    We prove \eqref{EE:1b} by repeating the arguments of Step~1.
      In view of \eqref{EE:1a} and~\eqref{EE:RET},  we have
\begin{equation}\label{EE:0b}
\lag  Q_1 \hat \xi,p_1(t)\rag=0,\quad \text{  $t\in J_\tau$}
\end{equation}  for any  $ \hat\xi=(0,\xi) \in \{0\}\times   G_2(\ty)$. 
  Taking the derivative   in time, we obtain 
 $$
\lag L_1( Q_1 \hat \xi) +b_1(\tilde v(t),Q_1 \hat \xi) + Q_1 \hat q(t) ,p_1(t)\rag =0, \quad\text{  $t\in J_\tau$},
$$ where 
$q(t)= L_2\xi+B_2(\tilde v(t),\xi)$ and $\hat q(t)=(Q_1 \hat \xi,q(t))\in H$,  $t\in J_\tau$.
 Taking another derivative and using the equality 
 $$\frac{\dd }{\dd t} Q_1 \hat q(t)=Q_1  (0,B_2(\dot{\tilde v}(t),\xi) ),$$ we get
\begin{align*}
\lag b_1(\dot {\tilde v}(t),Q_1 \hat \xi) + Q_1  (0,B_2(\dot{\tilde v}(t),\xi) ) ,p_1(t)\rag +\lag r(t) ,\dot p_1(t)\rag =0, \quad\text{  $t\in J_\tau$},
\end{align*}where $r(t)=L_1(Q_1 \hat \xi) +b_1(\tilde v(t),Q_1 \hat \xi)+Q_1 \hat q(t) $.
Now, we use the   equations for $\tilde v$ and $p_1$:
\begin{align*}
-\lag b_1(L_1\tilde v(t)+B_1(\tilde v(t))&+Q_1  \tilde u(t)-h_1 , Q_1 \hat \xi)) \\&- Q_1  (0,B_2(L_1\tilde v(t)+B_1(\tilde v(t))+Q_1  \tilde u(t)-h_1 ,\xi) ) ,p_1(t)\rag \\&+\lag   L_1 r(t)+b_1(\tilde v(t), r(t))  +Q_1\hat r (t)  ,p_1(t)\rag  =0, \quad\text{  $t\in J_\tau$},
\end{align*}where $\hat r (t)=( r(t) , L_2 q(t)+b_2(\tilde u(t), \hat q(t)))$. 
  Combining this with    \eqref{3.7}, \eqref{E:ffgdh}, and \eqref{ERTE},  we arrive at 
\begin{align*}
-\lag b_1(L_1\tilde v&+B_1(\tilde v)+Q_1  y-h_1 , Q_1 \hat \xi)) \\&- Q_1  (0,B_2(L_1\tilde v+B_1(\tilde v)+Q_1  y-h_1 ,\xi) ) ,p_1\rag \\&+\lag   L_1 r+b_1(\tilde v, r)   +Q_1 (r, L_2 q+B_2(\tilde v,  q)+B_2(Q_1\hat\xi, \pi_2 y)) ,p_1\rag  \\&
	 + \sum_{i=1}^d\lag - b_1( Q_1  \varphi_i , Q_1 \hat \xi)  +Q_1(0, \mathfrak{b}_2(\xi,\pi_2\varphi_i)) ,p_1\rag  \int_0^t\eta^i(s)\,\dd s=0.
\end{align*} Taking the derivative   in this equality and denoting
\begin{align*}
	\tilde a_i(t)&= \lag - b_1( Q_1  \varphi_i , Q_1 \hat \xi)  +Q_1(0, \mathfrak{b}_2(\xi,\pi_2\varphi_i)) ,p_1\rag  , \quad i\in[\![1,d]\!], \\
 	\tilde a_0(t)&=\frac{\dd}{\dd t} \Big(-\lag b_1(L_1\tilde v+B_1(\tilde v)+Q_1  y-h_1 , Q_1 \hat \xi))\\&\quad\quad -Q_1  (0,B_2(L_1\tilde v+B_1(\tilde v)+Q_1  y-h_1 ,\xi) ) ,p_1\rag \\&\quad\quad+\lag   L_1 r+b_1(\tilde v, r)   +Q_1 (r, L_2 q+B_2(\tilde v,  q)+B_2(Q_1\hat\xi, \pi_2 y)) ,p_1\rag  \Big)
	 \nonumber\\& \quad  +\sum_{i=1}^d\lag - b_1( Q_1  \varphi_i , Q_1 \hat \xi)  +Q_1(0, \mathfrak{b}_2(\xi,\pi_2\varphi_i)) ,\dot p_1\rag  \int_0^t\eta^i(s)\,\dd s , 
\end{align*}
we obtain 
$$
\tilde a_0(t)+\sum_{i=1}^d \tilde a_i(t)\eta^i(t)=0. 
$$ Again the functions~$\{\tilde a_i\}_{i=1}^d $ are continuously differentiable and $\tilde a_0$ is continuous, so the   observability of $\eta$ implies  that  $\tilde a_i(t) =0$ for $t\in J_\tau$ and~$i\in[\![0,d]\!]$. We have $\mathfrak{b}_2(\xi,\pi_2\varphi_i) \in \GG_2(\ty)$ for any $\xi\in  \GG_2(\ty)$. Hence, $Q_1(0, \mathfrak{b}_2(\xi,\pi_2\varphi_i))\in \GG_1(\ty)$, and  from~\eqref{EE:0b} it follows that
$$
\lag  Q_1(0, \mathfrak{b}_2(\xi,\pi_2\varphi_i)),p_1(t)\rag=0,\quad \text{  $t\in J_\tau$}.
$$Combining this with the equality  $\tilde a_i(t) =0$ for $t\in J_\tau$, we derive that  
$$
\lag b_1( Q_1  \varphi_i , Q_1 \hat \xi) ,p_1(t)\rag=0,\quad \text{  $t\in J_\tau$}.
$$Thus $\lag  \zeta ,p_1(t)\rag=0$ for any $\zeta\in \GG_1(\ty)$ and  $t\in J_\tau$.   The saturation assumption implies that \eqref{EE:1b} holds.

\section{Ergodicity of primitive equations}\label{S:mix} 

\subsection{Abstract result}

Here we formulate an abstract sufficient condition for  exponential mixing which is applied in the next section to the   randomly forced primitive~equations.    It is derived from Theorem~1.1 in~\cite{KNS-2018}.

%version of Theorem~1.1 in~\cite{KNS-2018}.  It is applied in the next section to establish exponential mixing for randomly forced primitive~equations.   

%In this section, we formulate a version of Theorem~1.1 in~\cite{KNS-2018}.   It is applied in the next section to establish exponential mixing for randomly forced primitive~equations.   

Let $H$ and $E$ be separable Hilbert spaces, let $\EE$ be a dense Banach subspace of~$E$, and let $X$ and $\KK\subset\EE$ be compact sets in $H$ and $E$, respectively. Assume that $S:X\times\KK\to X $ is a continuous mapping,  $\{\eta_k\}$ is a sequence of    i.i.d. random variables in $\EE$ with common law $\ell$ and   $\KK=\supp \ell$, and consider  a random sequence defined by
$$
u_k=S(u_{k-1},\eta_k),~k\ge1, \quad    u_0=u \in X.
$$  Then $(u_k, \pP_u)$, $u\in X$ is a Markov family in $X$, let~$\PPPP_k$ and~$\PPPP_k^*$ be the associated Markov operators. A measure~$\mu\in \PP(X)$ is said to be stationary for~$(u_k, \IP_u)$ if~$\PPPP_1^*\mu=\mu$. Recall that the dual-Lipschitz  metric  on~$\PP(X)$ is defined~by
$$
	 \|\mu_1-\mu_2\|_{L}^*=\sup_{f\in L(X), \, \|f\|_L\le1}\left|( f,\mu_1)-( f,\mu_2)\right|,
$$
where  $( f,\mu)=\int_Xf(u)\,\mu(\dd u)$ and  $L(X)$ is the space of   functions $f:X\to \R$ such~that 
$$
\|f\|_L=\sup_{u\in X} |f(u)| +\sup_{0<\|u-v\|_H\le 1}\frac{|f(u)-f(v)|}{\|u-v\|_H}<\infty.
$$
\begin{theorem}\label{T:crit}
Assume that the following   conditions hold.	
\begin{itemize}
\item[\hypertarget{H1}{\bf(H$_1$)}] 
{\sl There is a Banach space $V$ compactly embedded into  the space  $H$ such that~$X \subset   V$. There is an open set~$\OO=\OO_H\times \OO_E$ in $H\times E$ containing~$X\times \KK$  and 
 an extension $\tilde S:  \OO\to V$ of $S$ that is   twice continuously differentiable with    derivatives  that are bounded on bounded subsets of $  \OO$.~Moreover, for any~$u\in \OO_H$, the mapping $\eta\mapsto \tilde S(u,\eta)$, $\OO_E\to H$  is analytic, and all the derivatives $(D_\eta^j\tilde S)(u,\eta)$ are continuous in~$(u,\eta)$ and are bounded on bounded subsets of~$  \OO$.}

\item[\hypertarget{H2}{\bf(H$_2$)}] 
{\sl There are $a\in (0,1)$, $\hat \eta \in \KK $, and $\hat u \in X$ such that
$$
\|S(u,\hat \eta)-\hat u\|_H\le a\|u-\hat u\|_H \quad\text{for any $u\in X$.}
$$ }
 \item[\hypertarget{H3}{{\bf(H$_3$)}}]{\sl  For any $u\in X$ and~$\ell$-a.e. $\eta \in E$, the image of the linear mapping     $(D_\eta \tilde S)(u,\eta):E\to H$ is dense in~$H$. }

\item[\hypertarget{H4}{\bf(H$_4$)}]
{\sl The random variables $\eta_k$ are of the form
$\eta_k=\sum_{j=1}^\infty b_j\xi_{jk}e_j,$ where   $\{e_j\}$ is an orthonormal basis in $E$  such that $e_j\in \EE$ and 
$\sup_{j\ge1} \|e_j\|_\EE<\ty,$  $\{b_j\}$~are non-zero numbers  satisfying 
$\sum_{j=1}^\infty b_j^2 <\infty,$  and $\{\xi_{jk}\}$ are independent scalar random variables~with Lipschitz-continuous density $\rho_j$ such that $\supp \rho_j\subset [-1,1]$.}
\end{itemize}
Then the   family $(u_k,\IP_u), u\in X$ is exponentially mixing, i.e. it 
has a unique stationary measure $\mu\in\PP(X)$, and  there are numbers
  $C>0$ and~$c>0$ such~that 
\begin{equation}\label{4.1}
\|\PPPP_k^*\la-\mu\|_{L}^*
\le C  e^{-c k} , \quad    k\ge  0
\end{equation}for any initial measure  $\la\in \PP(X)$. \end{theorem}

 This theorem is a slight modification of Theorem~1.1 in~\cite{KNS-2018}.~The   difference is in Condition~(\hyperlink{H1}{H$_1$}) which is a localized version of the condition used in~\cite{KNS-2018}.~Indeed, it is not clear whether in the case of primitive equations this regularity condition holds with $\OO=H\times E$ (see Theorem~3 in~\cite{PMS}). Condition~(\hyperlink{H2}{H$_2$}) is usually satisfied   with $\hat u=0$  if the origin  is an exponentially stable
equilibrium for the unforced equation   and $0\in \KK$. Condition~(\hyperlink{H3}{H$_3$}) is a H\"ormander-type condition, and 
(\hyperlink{H4}{H$_4$})  is quite usual decomposability assumption.
We refer the reader to Section~1 in \cite{KNS-2018} for a detailed discussion of these conditions and for a short scheme of the proof of the original version of the theorem.

 \begin{proof}[Proof of Theorem \ref{T:crit}]
  Truncating   the mapping $\tilde S$, we easily obtain   an extension   
$\hat   S:H\times E\to V$ of~$S$   satisfying  (\hyperlink{H1}{H$_1$})   with   $\OO=H\times E$. Note that the family $(u_k,\IP_u), $ $u\in X$ does not change if we replace  $S$ by its extension $\hat  S$.   In view of Conditions~(\hyperlink{H1}{H$_1$})-(\hyperlink{H4}{H$_4$}),   the hypotheses of~Theorem~1.1 in~\cite{KNS-2018} are satisfied for the random dynamical system $u_k=\hat  S(u_{k-1},\eta_k)$. Applying that theorem, we     prove the mixing \eqref{4.1} for~$(u_k,\IP_u)$, $u\in X $.  
\end{proof}

\subsection{Application}

 In this section, we combine Theorems~\ref{T:3.3} and~\ref{T:crit}   to prove the exponential mixing for the randomly forced 3D primitive equations.~More precisely, we consider   Eq.~\eqref{1.3} with $h=0$ and  random   process  $\eta$     of the~form  
$$
	\eta(t)=\sum_{k=1}^\infty \I_{[k-1,k)}(t)\eta_k(t-k+1), \quad t\ge0,
$$
where  $\I_{[k-1,k)}$ is the indicator function of the   interval~$[k-1,k)$, $\{\eta_k\}$ is a sequence of    i.i.d.  random variables~in the space~$\EE=L^\ty(J,\HH)$, $J=[0,1]$, and~$\HH \subset U $ is a finite-dimensional subspace. In~what follows, we   denote by~$\ell$ the law of the random variable $\eta_k$ and assume that~$\KK=\supp \ell$ is compact in $\EE$. 
The~restriction to integer times   of the  solution  of Eq.~\eqref{1.3}   satisfies the~relation~$u_k=S_1(u_{k-1},\eta_k),$~$ k\ge1 	$ and
defines  a family of Markov processes~$(u_k, \IP_u)$ parametrised by the initial condition  $u_0= u \in V$.    
The following lemma is proved by   using   standard arguments based on  dissipative and  regularizing properties of PEs. 

%in~\cite{CT-2007, KZ-2007, J-2007, P-2006}; see  Theorem~3.1 in~\cite{Chu-14} for the case   when the PEs are perturbed by a bounded random kick force.

%{\red The following lemma is proved by   usual arguments based on dissipative and regularising properties of PEs established in~\cite{CT-2007, KZ-2007, J-2007, P-2006}; see  Theorem~3.1 in~\cite{Chu-14} for the case   when the PEs are perturbed by a bounded random kick force. }
\begin{lemma}\label{L:4.1}
	The family $(u_k, \IP_u)$ admits    a closed invariant absorbing set $X$  in~$U$  in the sense that, for   any $R>0$, there is an integer $k_0=k_0(R)\ge0$ such~that 
	\begin{align*}
			\pP_u\{u_k\in X,~k\ge k_0\}&=1 \quad \text{for $u\in B_V(0,R)$},\\
						\pP_u\{u_k\in X,~k\ge 0\}&=1 \quad \text{for $u\in X$}.
\end{align*}
\end{lemma} 
 %We shall consider    the restriction of $(u_k, \IP_u)$ to~$X$ which is compact in $V$.
 The following theorem is a more detailed version of
Theorem~\hyperlink{C}{C} formulated in the Introduction.
  \begin{theorem}  \label{T:4.1A}
	Let a finite-dimensional subspace $\HH\subset U$  be saturating in the sense of Definition~\ref{D:3.1},  and assume that  the following two conditions are fulfilled. 
\begin{description}
	\item [\bf Decomposability.] The random variables $\eta_k$ are of the form
$
\eta_k=\sum_{j=1}^\infty b_j\xi_{jk}e_j,
$
where   $\{e_j\}$ is an orthonormal basis in the Hilbert space~$E=L^2(J,\HH)$ such that 
$\sup_{j\ge1} \|e_j\|_{L^\ty(J,\HH)}<\ty,$  $\{b_j\}$ are non-zero   numbers  satisfying 
$\sum_{j=1}^\infty b_j^2 <\infty,$  and $\{\xi_{jk}\}$ are independent scalar random variables~with Lipschitz-continuous density $\rho_j$ such that $\supp \rho_j\subset [-1,1]$ and $\rho_j(0)>0$.
	\item [\bf Observability.]
There is $\tau\in(0,1)$ such that  the law  $\ell'$ of the restriction of the random variable  $\eta_k$ to the interval~$J_\tau$    is observable.\end{description}	
Then the   family $(u_k,\IP_u), u\in X$   is exponentially mixing.   
\end{theorem}
 \begin{proof} 
  By Theorem~3 in~\cite{PMS}, there is an open set~$\OO_E$ in $E$ containing~$\KK$   and  an extension $\tilde  S:H^2\times \OO_E\to H^3$ of $S_1$ that is    twice continuously differentiable  with  derivatives  that are bounded on bounded subsets of $H^2\times \OO_E$. Moreover,
  for any~$u\in H^2$, the mapping $\eta\mapsto \tilde S(u,\eta)$, $\OO_E \to H^2$ is analytic   and
   the derivatives~$(D^j_\eta \tilde S)(u, \eta)$~are continuous in~$(u, \eta)$ and   bounded on bounded subsets of~$H^2 \times \OO_E$. Thus~Condition~(\hyperlink{H1}{H$_1$}) in Theorem~\ref{T:crit} is verified with $H=H^2$,~$V=H^3$, and $\OO=H^2\times \OO_E$.

 Next,  for any    $\delta>0$, we define a norm     on $H^2$ by $|u|_\delta=\left(\|u\|^2+\delta\|u\|_2^2\right)^{1/2}$. Then  for~any bounded set~$B\subset H^2$, there are numbers $\delta>0$ and $a\in(0,1)$ such~that 
\begin{equation} \label{4.2}
	\bbar S_1(u)\bbar_\delta \le a\bbar u\bbar_\delta\quad\mbox{for $u \in B$}, 
\end{equation} where $S_1(u)=S_1(u,0)$.
This inequality with $B=X$  shows that    Condition~(\hyperlink{H2}{H$_2$})  is verified with $\hat u=0$,~$\hat \eta=0$, and the norm   $|\cdot|_\delta$.
To prove~\eqref{4.2},   we use the following inequalities     (see \cite{CT-2007, J-2007, P-2006, PMS}):
\begin{align*} 
\|S_1(u)\|&\le q\|u\|\quad\mbox{for $u\in H$},   \\ 
  \|S_1(u)\|_2&\le C_B\|u\|\quad\mbox{for $u\in B$},  	 
\end{align*}
where $q\in(0,1)$ and~$C_B>0$.   
These inequalities imply that 
$$ 
|S_1(u)|_\de^2 
=\|S_1(u)\|^2+\de\, \|S_1(u)\|_2^2
\le q^2\|u\|^2+\de\, C_B^2\|u\|^2
\le \left(q^2+\de C_B^2\right)\,|u|_\de^2. 
$$
Choosing $\delta>0$ so small that $a^2=q^2+\delta C_B^2<1$, we obtain~\eqref{4.2}.~Finally, Condition~(\hyperlink{H3}{H$_3$}) is established in Theorem~\ref{T:3.3}, and Condition~(\hyperlink{H4}{H$_4$}) is verified by the decomposability hypothesis. Applying Theorem~\ref{T:crit}, we complete the~proof.\end{proof}
 Theorem~\ref{T:4.1A} is formulated for initial measures $\lambda$ supported by $X$. As a   consequence of  Lemma~\ref{L:4.1} and Theorem~\ref{T:4.1A}, we obtain the following result. \begin{corollary}Under the conditions of Theorem~\ref{T:4.1A}, the measure $\mu$ is the unique   stationary measure for the family~$(u_k, \pP_u)$ in  $ \PP(V)$.	Moreover, inequality~\eqref{4.1}~holds for any~$R>0$,   $\la\in \PP(V)$ with~$\supp \la \subset B_V(0,R)$, and $k\ge k_0$. \end{corollary}

 \section{Saturating subspaces}
\label{S:4}

In this section, we show that the ten-dimensional subspace defined in the Introduction is  saturating in the  sense of both     Definitions~\ref{D:2.2} and~\ref{D:3.1}.~We also give an example of $V$-saturating subspace.

\subsection{$H$-saturating subspace}\label{S:4.1}

Let us consider the   subspace 
 \begin{equation}\label{5.1}
\HH= \text{span} \left\{(0,\phi_i): i\in  [\![1,10]\!]  \right\}\subset H,
\end{equation} where $\phi_i$  are the    eigenfunctions of the operator $L_2$ given by 
\begin{gather*}
\phi_1 =\cos  x\sin z,~\phi_2=\sin x\sin z,~\phi_3=\cos y\sin z,~\phi_4=\sin y \sin z ,\\  
\phi_5 =\sin z,~\phi_6=\cos  2x\sin z,~\phi_7=\sin 2x\sin z,~\phi_8=\cos 2y\sin z,\\\phi_9 =\sin 2y \sin z,~\phi_{10}=\sin 2z.
\end{gather*}
 \begin{theorem}\label{T:4.1}
	The subspace $\HH$ is   $H$-saturating in the sense of Definition~\ref{D:2.2}.
\end{theorem} 
To prove this theorem, we  introduce    the following two   orthogonal   bases:
\begin{description}
	\item [$\bullet$]     in $H_1$,   composed of eigenfunctions  of the operator   $L_1$:
	\begin{gather*}
	m\, c_m(x,y) \cos pz, \,\, m\, s_m(x,y) \cos pz, \,\, m^\bot c_m(x,y) \cos pz , \,\,   m^\bot s_m(x,y) \cos pz ,  \\ m^\bot c_m(x,y), \,\, m^\bot s_m(x,y), \,\, \imath\cos pz, \,\,\jmath\cos  pz  \quad m \in \Z^2_*,~p\ge1;
\end{gather*}
 	\item [$\bullet$]      in $H_2$,   composed of eigenfunctions  of the operator   $L_2$:
 	$$
c_m(x,y)\sin pz, \quad  s_m(x,y)\sin pz, \quad \sin pz  \quad m \in \Z^2_*,~p\ge1,
$$\end{description}	
  where 
   we denote $m^\bot=(-m_2,m_1)$,   $\imath=(1,0)$,   $\jmath =(0,1) \in \R^2$,  and
$$
c_m(x,y)= \cos(m_1x+m_2y), \quad\quad  s_m(x,y)= \sin(m_1x+m_2y).
$$The following two propositions are established in the next two subsections. 
\begin{proposition}\label{PP:5.2}
	Any vector of the     basis   in $H_1$ belongs\footnote{Here $\overline{\HH_1(\ty)}^{H_1}$ is the closure  of~$\HH_1(\ty)$ in $H_1$.} to $\overline{\HH_1(\ty).}^{H_1}$
\end{proposition}

\begin{proposition}\label{PP:5.3}
	Any vector of the     basis   in $H_2$ belongs to ${\GG}_2(\ty)$. 
\end{proposition}
  These   propositions  and the   inclusion $ \GG_2(\ty) \subset \HH_2(\ty)$  readily imply that~$\HH(\ty)$ is dense in    $H_1\times V_2$ and prove     Theorem~\ref{T:4.1}.

\subsubsection{Saturation in $\te$-component}\label{S:5.1.1}

In this subsection, we give a proof of Proposition~\ref{PP:5.3}.
\begin{proof}[Proof of Proposition~\ref{PP:5.3}] {\it Step~1.}~We first    show that   $\sin  x \sin pz\in  \GG_2(\ty)$  for   any  integer $p\ge1$.~Indeed, as~$\phi_2, \phi_5  \in   \GG_2$, we have  $\mathfrak{b}_2(\phi_2,\phi_5)\in {\GG}_2(1)$.  
   The equalities  
 \begin{gather*}
Q_1(0,\phi_2)= \imath \cos x\cos z, \quad Q_1(0,\phi_5)= 0, \\
	\mathfrak{b}_2( \phi_2 , \phi_5)=B_2(Q_1(0,\phi_2), \phi_5)=	\frac{1}{2}\sin x\sin 2z
\end{gather*}  
imply that
   $ \sin x\sin 2z\in {\GG}_2(1)$.   A simple computation shows that
 \begin{align*}
 Q_1(0, \sin x\sin 2z)&=\frac12\imath \cos x\cos 2z,\\
	\mathfrak{b}_2(\sin x\sin 2z,\phi_5)& =B_2(Q_1(0, \sin x\sin 2z),\phi_5)\\&= \frac14\sin x\sin 2z \cos z\\&=\frac1{8}\left(\sin x\sin z + \sin x\sin 3z  \right)\in {\GG}_2(2).
\end{align*}
 This implies that   $\sin x\sin 3z \in {\GG}_2(2)$.~Iterating this argument, we see that   $ \sin x\sin pz \in {\GG}_2(\ty)$ for any $p\ge1$.~In a similar way, we can prove that      $\cos x \sin pz$, $\cos y \sin pz$,   $\sin y \sin pz\in {\GG}_2(\ty)$ for any $p\ge1$.

Repeating the above arguments and using the fact that $\sin2x\sin z$, $\cos2x\sin z$, $\sin2y\sin z$, $\cos2y\sin z\in {\GG}_2$, we can obtain also that $\sin2x\sin pz$, $\cos2x\sin pz$, $\sin2y\sin pz$, $\cos2y\sin pz\in {\GG}_2(\ty)$ for any $p\ge 1$.
 
{\it Step~2.}  
 	  Let us show   that $\sin pz \in {\GG}_2(\ty)$    for  any  $p\ge1$.  The equalities  
 	\begin{align*}
 		 \imath\cos x \cos nz&=Q_1(0, n \sin x \sin nz), \\
 		   \imath\sin x \cos nz&=Q_1(0, -n\cos x\sin nz)  
 	\end{align*} and the fact that $\cos x\sin nz$, $\sin x \sin nz\in {\GG}_2(\ty) $
 	imply that 
 	   	 \begin{align*}
   		\bB_2(n\cos x\sin nz,&\phi_1)+\bB_2(n \sin x \sin nz,\phi_2)\\&=\frac{n^2-1}{2n} \left((n-1)\sin(n-1)z-(n+1)\sin(n+1)z\right)\in {\GG}_2(\ty).
   \end{align*}
   Thus $n \sin nz-p\sin pz \in {\GG}_2(\ty)$ for any  integers   $n>p\geq1$ that are both even or both odd. 
   As   $\sin z$, $\sin 2z\in {\GG}_2$, we obtain~$\sin pz \in {\GG}_2(\ty)$ for any~$p\ge1$.

 {\it Step~3.}  
   Now we prove    the following property $P(m)$   for any $m\ge1$:
   \begin{description}
   \item[$P(m): $] for any   $p\ge1$ and $q\in[\![1,m]\!]$, we have 	 $ \sin qx \sin pz \in {\GG}_2(\ty)$. 
   \end{description}
    We argue by induction on $m$. The cases $m=1,2$ are considered in Step~1. Assuming  that   $P(m)$ is true for $m\ge 2$, we  prove it for $m+1$. Note~that 
      \begin{align}
   \imath\cos mx \cos nz&=\frac nm Q_1(0,  \sin m x \sin nz), \nonumber \\
      \bB_2\left(\frac n{m} \sin m x \sin nz,\phi_1\right) &= 
      \frac{(m-n)(m+n^2)}{4mn}\sin(m+1)x\sin (n+1)z 
     \nonumber \\&\quad \!+\!\frac{(m+n)(m+n^2)}{4mn}\sin(m+1)x\sin (n-1)z
      \nonumber\\&\quad \!+\!\frac{(m+n)(m-n^2)}{4mn}\sin(m-1)x\sin (n+1)z
      \nonumber\\&\quad \!+\!\frac{(m-n)(m-n^2)}{4mn}\sin(m-1)x\sin (n-1)z  \in    {\GG}_2(\ty) 	\label{5.2}
   \end{align}for any $n\ge1$. By the induction hypothesis, we have $\sin(m-1)x\sin (n+1)z$, 
      $\sin(m-1)x\sin (n-1)z \in   {\GG}_2(\ty) $. Thus  \eqref{5.2} implies that 
\begin{equation}\label{5.3}
 (m-n) \sin(m+1)x\sin (n+1)z 
      + (m+n) \sin(m+1)x\sin (n-1)z
       \in   {\GG}_2(\ty)
\end{equation}Taking here  $n=1$, we obtain $\sin(m+1)x\sin 2z
       \in   {\GG}_2(\ty).$ It follows~that 
     \begin{align*}
       \bB_2\bigg(\frac2{m+1} \sin(m+1)x\sin 2z,\phi_5\bigg)\!    =\!\frac{(m+1)}{4}\sin (m+1)x(\sin3z+\sin z)\!\in\! {\GG}_2(\ty).
      \end{align*}
       Thus
       \begin{equation}\label{5.4}
       	\sin (m+1)x \sin3z+ \sin (m+1)x \sin z \in  {\GG}_2(\ty).
       \end{equation}
Taking $n=2$ in \eqref{5.3}, we get 
$$ (m-2) \sin(m+1)x\sin 3z 
      + (m+2) \sin(m+1)x\sin  z
       \in   {\GG}_2(\ty).
$$
This and \eqref{5.4} imply that $\sin(m+1)x\sin  z
       \in   {\GG}_2(\ty)$.
       Repeating the argument of Step~1, we show  that $\sin(m+1)x\sin  pz
       \in   {\GG}_2(\ty)$ for any~$p\ge1$. Thus $P(m)$ is true for any $m\ge1$

   %To show that $\cos 2x \sin pz \in \hH_2$, we use the equality
  % \begin{align*}
  % 	B_2(v_n^c,\phi_2)+	B_2(v_n^s,\phi_1)=\frac{1}{2n}\cos2x \left((n-1)\sin(n+1)z-(n+1)\sin(n-1)z\right) .
  % \end{align*}
% As $v_n^c, v_n^s \in \HH_1(\ty)$,    we have 
% $$
% \frac1{n+1}\cos2x\sin(n+1)z- \frac1{n-1}\cos2x\sin(n-1)z \in \HH_2(\ty), \quad n\ge1.
% $$This implies that 
 % $$
 %\frac1{n}\cos2x\sin nz- \frac1{p}\cos2x\sin p z \in \HH_2(\ty), \quad n>p\ge1.
% $$
%Thus  properties  (i) and (ii) are satisfied for   $\cos 2x \sin pz$  with 
%$$
%g^{c}_{n,(2,0),p}= - \frac {p}{n}\cos2x\sin nz,
% $$so~$\cos 2x \sin pz\in \hH_2$ for any $p\ge1$.

In a similar way,      $\cos mx \sin pz$, $\cos my \sin pz$,   $\sin my \sin pz$ belong to ${\GG}_2(\ty)$ for any $m,p\ge1$.

{\it Step~4.}  In this step, we show that $  s_m(x,y)\sin pz\in {\GG}_2(\ty)$  for any $p\ge1$ and~$m=(m_1,m_2)\in \Z_*^2$.  We confine ourselves to the case~$m_2\ge0$,  the   case~$m_2<0$ being similar. Arguing by induction on $m_2$,  we prove   the following property:
   \begin{description}
   \item[$P'(m_2):$] for any   $p\ge1$, $m_1\in \Z $,   and $q\in [\![0,m_2]\!]$, we have
    $ s_m(x,y)\sin pz \in {\GG}_2(\ty),$ where $m=(m_1, q)$.   
   \end{description}The case $m_2=0$ is considered in Step~3.   Assuming that $P'(m_2)$ is true for~$m_2\ge 0$, let us prove  it for $m_2+1$. 
   We first consider the case $m_1\neq \pm 1$. 
    Let     \begin{align*}
   \theta_1&=-n\cos(m_1x+m_2y)\sin nz,\\
      \theta_2&=n\sin(m_1x+m_2y)\sin nz.		
   \end{align*}Using the equalities   
   \begin{align*}
   Q_1(0,\te_1)&=m s_m(x,y) \cos nz,	\\
      Q_1(0,\te_2)&=m c_m(x,y) \cos nz,	
   \end{align*}we get  
   \begin{align}
   \bB_2(\te_1,\sin y\sin z)&-	 \bB_2(\te_2,\cos y\sin z)= a_1\sin(m_1x+(m_2+1)y)\sin(n+1)z\nonumber\\&\quad +a_2\sin(m_1x+(m_2+1)y)\sin(n-1)z\in {\GG}_2(\ty), \label{5.5}
   \end{align}where 
   \begin{align*}
   a_1&=\frac1{2n}\left(n^3-n(n-1)m_2-m_1^2-m_2^2\right),	\\
   a_2&=-\frac1{2n}\left(n^3+n(n+1)m_2+m_1^2+m_2^2\right).
   \end{align*}
   Taking $n=1$ in \eqref{5.5}, we obtain 
   $$
   (1-m_1^2-m_2^2)\sin(m_1x+(m_2+1)y)\sin2z\in {\GG}_2(\ty).
   $$As $m_1\neq \pm 1$, we have $1-m_1^2-m_2^2\neq 0$, so 
$
 \sin(m_1x+(m_2+1)y)\sin2z\in {\GG}_2(\ty).
$ The  latter implies that
     \begin{align}
       \bB_2(  2\sin(m_1x+(m_2+1)y)\sin2z,\phi_5)    =&\frac{m_1^2+(m_2+1)^2}{4}\sin(m_1x+(m_2+1)y)\nonumber\\&\quad \times (\sin3z+\sin z)\in {\GG}_2(\ty). \label{5.6}
      \end{align}On the other hand, taking $n=2$ in \eqref{5.5}, we get 
\begin{equation}\label{5.7}
a_1\sin(m_1x+(m_2+1)y)\sin3z  +a_2\sin(m_1x+(m_2+1)y)\sin z\in {\GG}_2(\ty).
\end{equation}When $n=2$, we have
$
a_1-a_2=4+m_2\neq 0, 
$ since $m_2\ge 0$. Combining~\eqref{5.6} and~\eqref{5.7}, we see that 
$ \sin(m_1x+(m_2+1)y)\sin z\in {\GG}_2(\ty).
$ Now applying the argument of Step~1, we infer  that  $ \sin(m_1x+(m_2+1)y)\sin pz\in {\GG}_2(\ty)$ for any~$p\ge1$ and $m_1\neq \pm1$.  

Finally, computing the term $\bB_2(\sin(\pm2x+(m_2+1)y)\sin nz,\phi_1)$, one easily shows that $ \sin(\pm x+(m_2+1)y)\sin pz\in {\GG}_2(\ty)$ for any $p\ge1$.
  Thus $P'(m_2)$ holds for any $m_2\ge0$, and
we conclude that that $  s_m(x,y)\sin pz\in {\GG}_2(\ty)$  for any~$p\ge1$. The proof of $  c_m(x,y)\sin pz \in {\GG}_2(\ty)$ is similar.

Thus all the vectors of the basis in $H_2$     belong to ${\GG}_2(\ty)$. This completes the proof of Proposition~\ref{PP:5.3}. 
\end{proof}

\subsubsection{Saturation in $v$-component}\label{S:5.1.2}

      Here we prove Proposition~\ref{PP:5.2}.  We first establish the following lemma. 
        	 \begin{lemma}\label{C:4.3}For any $m \in \Z^2_*$, $p\ge1$, and  $i\in [\![1,4]\!]$, the following    properties hold:
      	\begin{enumerate}
      	\item[(a)]  the functions $m \,c_m(x,y)\cos pz$, $m \, s_m(x,y)\cos pz$  belong to $ \HH_1(\ty)$; 
                \item[(b)] the functions 
                	$$\psi_{m,p,i}^c=b_1(m \,c_m(x,y)\cos pz,\psi_i), \,\, \psi_{m,p,i}^s=b_1(m \,s_m(x,y)\cos pz,\psi_i)$$
                belong to $ \HH_1(\ty)$, where   
                  	$$\psi_1=\iota  \cos  x \cos z,\,\, \psi_2=\iota \sin  x \cos z,\,\,    \psi_3=\jmath \cos  y \cos z, \,\, \psi_4=\jmath \sin  y \cos z.$$

      	\end{enumerate}
  	\end{lemma}
  \begin{proof}   By Proposition~\ref{PP:5.3},  we know that   $s_m(x,y)\sin pz$, $c_m(x,y)\sin pz \in \HH_2(\ty)$. Recall that 
  $ Q_1(0,\te)\in \HH_1(\ty)$ for any $\te \in \HH_2(\ty)$. So property (a) follows from the~equalities
  \begin{align*}
  	 	m \,c_m(x,y)\cos pz &=   Q_1(0,p s_m(x,y)\sin pz), \\
  	 	 	m \,s_m(x,y)\cos pz &=  Q_1(0,-p c_m(x,y)\sin pz) , \quad  m \in \Z^2_*,\,p\ge1.  \end{align*} 
 To prove (b), we take any~$\e>0$ and   $a\in \R$  and  
  note that
 \begin{align*}
	 B_1(\e \,  m \,c_m(x,y)\cos pz + a\e^{-1}\psi_i)&= \e^{-2}B_1(a\psi_i)  + \e^2B_1( m \,c_m(x,y)\cos pz)\nonumber\\ &\quad +	b_1( m \,c_m(x,y)\cos pz, a\psi_i). 
\end{align*} 
 Using the fact that $B_1(\psi_i)=0$,  we obtain    the following  limit   in $V_1$ as~$\e\to 0^+$:
 $$
 B_1(\e  m \,c_m(x,y)\cos pz + a\e^{-1}\psi_i) \to b_1( m \,c_m(x,y)\cos pz, a\psi_i)=a \psi^c_{m,i}.
 $$  As $a\in \R $ is arbitrary, this shows that~$ \psi_{m,p,i}^c \in \HH_1(\ty)$.  Repeating these arguments with $s_m(x,y)$ instead of $c_m(x,y)$, we prove that   $ \psi_{m,p,i}^s \in \HH_1(\ty)$.  
 	\end{proof}

\begin{proof}[Proof of Proposition~\ref{PP:5.2}]
 	{\it Step~1.}~Let us    show that   
	  $\imath \cos pz  \in  \overline{\HH_1(\ty)}^{H_1}$ for any~$p\ge1$. To this end,
  we take any $n\ge2$ and compute the term  
 \begin{align}
	 b_1(\imath\cos x\cos nz, \psi_2)&=   \frac{1}{2}    \imath \cos 2x\left(\cos (n+1)z+\cos (n-1)z\right)\nonumber\\&\quad+ \frac{1+n^{2}}{2n}  \imath\cos 2x \left(\cos (n-1)z-\cos (n+1)z\right)\nonumber \\&\quad +\frac{n^{2}-1}{2n} \imath\left(\cos (n-1)z-\cos (n+1)z\right).\nonumber
\end{align}
By property (a) in Lemma~\ref{C:4.3}, we have    $ \imath \cos 2x\cos (n\pm 1)z\in \HH_1(\ty), $ and by property (b), that $ b_1(\imath\cos x\cos nz, \psi_2)\in \HH_1(\ty). $ It follows that 
$$
     \imath\left(\cos (n-1)z-\cos (n+1)z\right)\in \HH_1(\ty) \quad \text{for any $n\ge2$},
$$
     so 
           \begin{equation}\label{5.8}
         \imath\left(\cos pz-\cos qz\right)\in \HH_1(\ty) ,  
 \end{equation}
  provided that $p,q\ge1$ are both odd or both even. Passing to the limit as $q\to \ty$, we see that, for any $p\ge1$, the function  $\imath \cos pz$ is in the   $L^2$-weak closure of $\HH_1(\ty)$, hence in   $ \overline{\HH_1(\ty) }^{H_1}$, since $\HH_1(\ty)$ is a vector space.

Computing the term  
$
b_1(\jmath\cos y\cos pz,  \psi_4)
$ and repeating the above arguments,    
we infer that   
	  $\jmath \cos pz \in   \overline{\HH_1(\ty).}^{H_1}$  
	   
	{\it Step~2.} In this step, we show that $m^\bot s_m(x,y)\cos pz$  belongs to  $\overline{\HH_1(\ty)}^{H_1}$   for any $m \in \Z^2_*$ and $p\ge0$.    	
		 Let us take any  $m=(m_1,m_2)\in \Z^2_*$ and $n\ge1$, and use the equality   
		\begin{align}
			b_1(s_{m+\jmath}(x,y)&\cos nz,   \psi_4 )+b_1(c_{m+\jmath}(x,y)\cos nz,  \psi_3 )\nonumber
\\&= -(m_2
+1) \Pi \left(s_m(x,y) \cos nz \cos z\right)\nonumber\\&\quad\, + \Pi \left(A_\jmath(m,n) \sin(m_1x+m_2y) \sin nz\sin z\right),\label{5.9}
\end{align}where $ A_\jmath(m,n) =n^{-1}( m_1 n^2, (m_2+1)n^2-m_1^2-(m_2+1)^2)$.~By   Lemma~\ref{C:4.3},  we have that the functions
  $ \Pi \left(s_m(x,y) \cos nz \cos z\right)$, $b_1(s_{m+\jmath}(x,y)\cos nz,   \psi_4 ),$ and $b_1(c_{m+\jmath}(x,y)\cos nz,  \psi_3 )$ belong  to $\HH_1(\ty).$ Hence, 
  \begin{equation}\label{5.10}
\Pi \left(A_\jmath(m,n) \sin(m_1x+m_2y) \sin nz\sin z\right) \in \HH_1(\ty).
\end{equation}       
		The vector $A_\jmath(m,n)$ is parallel to $m$ if and only if one of the following two conditions hold:
		\begin{enumerate}
		\item[$\bullet$]	$m_1=0$ and $m_2\neq 0$;
			\item[$\bullet$]	$m_1\neq 0$ and  $n^2=m_1^2+(m_2+1)^2.$
		\end{enumerate}
 Let us denote by $\aA_\jmath$ the   set of couples $(m,n)$ such that $A_\jmath(m,n)$ is non-parallel to $m$. From \eqref{5.10} we derive that 
		\begin{equation}\label{5.11}
		m^\bot s_m(x,y) (\cos (n+1)z-\cos (n-1)z)	 \in  \HH_1(\ty) 
		\end{equation} for any $(m,n)\in \aA_\jmath$.
 		 In a similar way, we compute   the sum
			\begin{align}
			b_1(s_{m+\iota}(x,y)&\cos nz,   \psi_2 )+b_1(c_{m+\iota}(x,y)\cos nz,  \psi_1)
\nonumber\\&= -(m_1+1) \Pi\left(s_m(x,y) \cos nz \cos z\right)\nonumber\\&\quad\, + \Pi \left(A_\iota(m,n) \sin(m_1x+m_2y) \sin nz\sin z\right),\label{5.12}
\end{align}where $A_\iota(m,n) =n^{-1}(  (m_1+1)n^2-m_2^2-(m_1+1)^2, m_2 n^2)$. As above, $A_\iota(m,n)$ is parallel to $m$ if and only if one of the following  conditions hold:
		\begin{enumerate}
		\item[$\bullet$]	$m_2=0$ and $m_1\neq 0$;
			\item[$\bullet$]	$m_2\neq 0$ and  $n^2=(m_1+1)^2+m_2^2.$
		\end{enumerate}Let    $\aA_\iota$ be the   set of   $(m,n)$ such that $A_\iota(m,n)$ is non-parallel to $m$.  From~\eqref{5.12} it follows that \eqref{5.11} holds for any $(m,n)\in \aA_\iota$.

 Let us go  back to \eqref{5.9}, and  replace $m$ by $-m$. We see that
$$
\Pi \left(A_\jmath(-m,n) \sin(m_1x+m_2y) \sin nz\sin z\right) \in \HH_1(\ty),
$$          and $A_\jmath(-m,n)$ is parallel to $m$ if one of the following  conditions hold:
\begin{enumerate}
		\item[$\bullet$]	$m_1=0$ and $m_2\neq 0$;
			\item[$\bullet$]	$m_1\neq 0$ and  $n^2=m_1^2+(m_2-1)^2.$
		\end{enumerate}We denote by $\aA_\jmath^-$ the   set of   $(m,n)$ such that $A_\jmath(m,n)$ is non-parallel to $m$. Again \eqref{5.11} holds for any $(m,n)\in \aA_\jmath^-.$ The set $\aA_\iota^-$ is defined in a similar way, by replacing $m$ by $-m$ in \eqref{5.12}. Then  \eqref{5.11} holds for any $(m,n)\in \aA_\iota^-.$
		
		It is easy to see that the union of the sets $\aA_\jmath^\pm$, $\aA_\iota^\pm$ is $\Z^2_*\times \N_*$, so   \eqref{5.11} holds for {\it any} $m \in \Z^2_*$ and $n\ge1$.  Iterating \eqref{5.11}, we obtain  
\begin{equation}\label{5.13}
m^\bot s_m(x,y) (\cos pz-\cos qz)	 \in  \HH_1(\ty),
		\end{equation} provided that the integers $p,q\ge0$  are   both even or   odd. Passing to the limit as $q\to \ty$, we conclude  that 
$
m^\bot s_m(x,y)  \cos pz 	 
$ belongs to the   $L^2$-weak closure of~$\HH_1(\ty)$, hence to   $ \overline{\HH_1(\ty) }^{H_1}$  for any~$m\in \Z^2_*$   and  $p\ge0$.

	 A similar argument shows that    $m^\bot c_m(x,y)\cos pz\in \overline{\HH_1(\ty)}^{H_1}$   for any $m \in \Z^2_*$ and $p\ge0$.   This completes the proof of Proposition~\ref{PP:5.2}.  
		\end{proof}

\subsection{$V$-saturating subspace}\label{S:4.3}

 We do not know whether the subspace $\HH$ defined in \eqref{5.1} is $V$-saturating. The main difficulty comes from Proposition~\ref{PP:5.2}, which   only implies   that  $\HH_1(\ty)$ is
dense in $H_1$. In this section, we add new vectors to $\HH$ with non-zero $v$-components. This results in a larger    $\HH_1(\ty)$-space  that contains    all the vectors of the basis in $H_1$. As a consequence, we get $V$-saturation property.  More precisely,  we  define the space 
$$
\tilde \HH=\text{span}\{(\tilde \phi_i,0), (0,\phi_j): i=1,  \ldots,  6,\,\, j=1, \ldots, 10\}\subset H,$$
 where the functions  $\phi_j$ are as in Section~\ref{S:4.1} and 
\begin{gather*}
\tilde \phi_1=\jmath  \cos z, \, \, \tilde \phi_2=\jmath  \cos 2 z, \,\,\tilde \phi_3=\imath   \cos z,\,\,  \tilde \phi_4 =\imath    \cos 2z,\\\tilde \phi_5=\jmath  \cos x, \, \, \tilde \phi_6=\jmath  \sin x.
\end{gather*}
\begin{theorem}\label{T:4.4}
	The subspace $\tilde \HH$ is   $V$-saturating in the sense of Definition~\ref{D:2.2}.
\end{theorem}
\begin{proof}Let  $\tilde \HH(j)$ be the subspaces defined by \eqref{2.4} and \eqref{2.5} with $\HH=\tilde \HH$, and let~$\tilde \HH_i(j)=\pi_i \tilde \HH(j)$, $i=1,2$.
From Proposition~\ref{PP:5.3} it follows that  $\tilde \HH_2(\ty)$ is~dense in~$V_2$.  The proposition will be proved if we show that any vector of the     basis   in $H_1$ belongs  to~$ \tilde \HH_1(\ty)$   (cf. Proposition~\ref{PP:5.2}). By Lemma~\ref{C:4.3}, we have 
   $m \,c_m(x,y)\cos pz$,  $m \, s_m(x,y)\cos pz\in  \tilde \HH_1(\ty)$ for any~$m\in \Z^2_*$ and~$p\ge1$.
Combining~\eqref{5.8}, the version of~\eqref{5.8}   with $\jmath$ instead of~$\iota$, and~the~assumption that~$(\tilde \phi_i,0) \in \tilde \HH$, $i=1,\ldots,4$, we obtain that~$  \imath\cos pz$, $\jmath\cos  pz \in  \tilde \HH_1(\ty)$  for any~$p\ge1$. For any $ m\in \Z^2_*$ and $p\ge1$, the following
  equality holds:
\begin{gather}
			b_1(s_{m+\iota}(x,y)\cos pz,   \tilde\phi_6)+b_1(c_{m+\iota}(x,y)\cos pz,   \tilde\phi_5)\nonumber
\nonumber\\ =   \Pi \left(A(m) \sin(m_1x+m_2y) \cos pz\right), \label{5.14}
\end{gather} where $A(m)=(-(m_1+1)m_2, m_1+1-m_2^2).$ Note that the vector  $A(m)$ is parallel to~$m$  if and only if $A(m)=0$, i.e., $m=(-1,0).$ Assume that $m\neq (-1,0)$. Since $B_1(\tilde\phi_5)=B_1(\tilde\phi_6)=0$, as in the proof of (b) in Lemma~\ref{C:4.3}, we   show that~$b_1(s_{m+\iota}(x,y)\cos pz$,   $\tilde\phi_6), b_1(c_{m+\iota}(x,y)\cos pz,   \tilde\phi_5) \in \tilde\HH_1(\ty)$. As $A(m)$ is non-parallel to $m$, from \eqref{5.14} we derive that $m^\bot s_m(x,y) \cos pz\in \tilde \HH_1(\ty)$ for any $p\ge1$. From \eqref{5.13} it follows   also that  $m^\bot s_m(x,y) \in \tilde \HH_1(\ty)$. Finally, if~$m=(-1,0)$, then $m^\bot s_m(x,y)=\tilde \phi_6(x)\in \tilde \HH $,   and   (cf. (b) in Lemma~\ref{C:4.3})
$$
b_1((1,1) \cos pz, \tilde \phi_5(x))=-m^\bot s_m(x,y) \cos pz\in \tilde \HH_1(\ty), \quad p\ge1. 
$$  
With similar arguments one proves also that  $m^\bot c_m(x,y)\cos pz \in \tilde \HH_1(\ty)$ for any~$m\in \Z^2_*$ and $p\ge0$. Thus any vector of the basis   in $H_1$ belongs  to~$ \tilde \HH_1(\ty)$. We conclude that  $ \tilde \HH_1(\ty)$ is dense in $V_1$ and $\tilde \HH(\ty)$ is dense in $V$.
   \end{proof}

 \subsection{Saturation   for   linearized system} 

Now we turn to the saturation property for the  linearized system.
 \begin{theorem}\label{T:5.6}
	The subspace $\HH$        defined by \eqref{5.1}  is saturating for     linearized system~\eqref{1.11}
	  in the sense of Definition~\ref{D:3.1}.
\end{theorem} 
\begin{proof}This theorem follows from the   proof of Theorem~\ref{T:4.1}.  Indeed, by Proposition~\ref{PP:5.3},  $\GG_2(\ty)$ is dense in $H_2$.
  The computations in Section~\ref{S:5.1.2} show that any vector   of the basis   in $H_1$ belongs  to~$\overline{ \GG_{1}(\ty).}^{H_1}$ We conclude that~$\GG(\ty)$ is dense in~$H$, so  $\HH$ is saturating in the sense of Definition~\ref{D:3.1}.
\end{proof}

\section{Proof of Proposition \ref{P:1.2}}
\label{S:6}

We confine ourselves to the proof  of limit \eqref{1.10}, which is relatively more complicated; see    Remark~\ref{R:R}.
  Let us take any~$u_0=(v_0,\te_0)\in H^4$  and   $\xi \in  H^5$ such that $\pi_1\xi=0$  and consider the  function
$
	w(t)= u(\delta t)-q(t),    
$
 where $u(t)=\RR_t(u_0,    \delta^{-1} \xi,0)$, $ q(t)=(q_1(t),q_2(t)),$  
\begin{align}
   q_1(t)&=v_0-t Q_1 \xi, \label{6.1} \\
    q_2(t)&=\te_0-t \left(L_2\zeta+B_2(v_0,\zeta)\right)+\frac {t^2}2 B_2(Q_1\xi,\zeta),\label{6.2}
\end{align}
and $\zeta=\pi_2 \xi$. 
   For any      $r>0$, we show that
   \begin{equation}\label{dddsss}
     w(1)\to 0 \quad \textup{in~$V$ as~$\delta\to 0^+$}
     \end{equation}
    uniformly with respect to
$u_0$ and $\zeta$  satisfying 
\begin{equation}\label{6.3}
 \|u_0\|_4+\|\zeta\|_5\le r.
\end{equation} Note that          $v=\pi_1 w $  is a solution of the following equation:  
 \begin{gather*}
  \p_t v   + \delta L_1(v+q_1)+ \delta\lag v+q_1, \nabla\rag (v+q_1)    - \delta\int_0^z \diver  (v +q_1)\,\dd \z \,
\p_z (v+q_1)  \nonumber\\+   \delta f   (v+q_1)^\bot  +
 \delta \nabla p_s -\delta \int_0^z \nabla(\te +q_2+\de^{-1}\zeta)\,\dd \z+ \p_t q_1  = \de h_1.   
 \end{gather*}  Using \eqref{6.1} and \eqref{E:qban}, we see that this equation is equivalent to
  \begin{gather}
  \p_t v   + \delta L_1(v+q_1)+ \delta\lag v+q_1, \nabla\rag (v+q_1)    - \delta\int_0^z \diver  (v +q_1)\,\dd \z \,
\p_z (v+q_1)  \nonumber\\+   \delta f   (v+q_1)^\bot  +
 \delta \nabla p_s -\delta \int_0^z \nabla(\te +q_2)\,\dd \z  = \de h_1.  \label{6.4} 
 \end{gather} In a similar way, $\te=\pi_2 w$ is a solution of the equation
 \begin{gather*}
  \p_t\te  +  \delta L_2 (\te +q_2 +\delta^{-1}\zeta)    +  \delta \lag v +q_1 ,  \nabla\rag (\te + q_2 +\delta^{-1} \zeta)  \nonumber   \\-  \delta\int_0^z \diver (v+q_1) \,\dd \z\,
\p_z (\te+q_2+ \delta^{-1} \zeta)  
+  \p_t q_2    =  \de h_2 ,   
\end{gather*} which, in view of
 \eqref{6.2} and \eqref{E:B2},  can be  rewritten         as   follows:
\begin{gather}
  \p_t \te  +  \delta L_2 (\te +q_2)    +  \delta \lag v +q_1 ,  \nabla\rag (\te + q_2) \nonumber   \\-  \delta\int_0^z \diver (v+q_1) \,\dd \z\,
\p_z (\te+q_2)+   \lag v ,  \nabla\rag \zeta  - \int_0^z \diver  v  \,\dd \z\,
\p_z   \zeta   
  =  \de h_2 .  \label{6.5}
\end{gather}
   The initial conditions are $ v(0) =0$ and  $ \te(0)=0$.   Using   Eqs.~\eqref{6.4} and~\eqref{6.5},   we   estimate the norms of $v$ and $\te$.  We begin with estimates for $L^2$-norms of $v$ and $\theta$, where we use standard energy methods. Then, 
   following the ideas of Cao and  Titi~\cite{CT-2007}, we  
  estimate the~$L^6$-norm of $v$ and $H^1$-norms of $v$ and $\te$   using specific  barotropic-baroclinic formulation of PEs. The latter consists in representing  the velocity field $v$ as follows $v=\bar v +\tilde v,$ where $\bar v$ is the vertical average of $v$ and $\tilde v$ is  the remaining part. The average $\bar v$ is two-dimensional, so it is treated by usual 2D NS methods, and  the advantage of the reminder $\tilde v$ is that it satisfies an equation without pressure.
    We provide a very detailed proof of limit~\eqref{dddsss} that is divided into nine steps.

 {\it Step~1.~$L^2$-estimate for $v$.}~The goal of the first two steps is to show the limit~$\|v(1)\|+\|\te(1)\| \to 0$  as~$\delta\to 0^+$. We start by taking  the scalar product in~$L^2$ of Eq.~\eqref{6.4} with~$v$ and integrating by parts:
  \begin{align}\label{6.6}
 \frac{1}{2}  \frac{\dd }{\dd t} \|  v\|^2+&\delta \nu_1\|\nabla  v\|^2+\delta\mu_1\|\p_{z} v\|^2 =  - \delta\left\lag L_1  q_1 , v\right\rag  -\delta\left\lag \lag v+q_1, \nabla\rag (v+q_1),  v \right \rag   \nonumber\\
&+ \delta\left\lag\int_0^z \diver  (v +q_1)\,\dd \z \,
\p_z (v+q_1),  v \right\rag   -  \delta \left\lag f   (v+q_1)^\bot,  v\right\rag   \nonumber\\&  -
 \delta\left\lag\nabla p_s,  v\right\rag +\delta\left\lag \int_0^z \nabla(\te +q_2)\,\dd \z ,  v \right\rag + \delta \left\lag h_1 ,  v \right\rag = \sum_{i=1}^7I_i.
\end{align} To estimate the terms   $I_1, I_4, I_6, I_7$, we integrate by parts and use\footnote{These   inequalities and assumptions   are used in almost all the    estimates below, so we will not mention them every time. The same letter $C$ is used to denote   constants which may change from line to line.}  the Cauchy--Schwarz and Young inequalities  and  the assumption that       $u_0
$  and $\zeta$  satisfy \eqref{6.3} and~$t\in[0,1]:$ \begin{align*}
|I_1|&\le \delta \|L_1 q_1 \| \, \|  v\|\le  C \delta      \|    v\|, \\
	| I_4|&= \delta |\lag f   q_1^\bot,   v \rag|\le     \delta \|fq_1\|\, \|  v\|\le   C\delta    \|  v\|,   \\
	|I_6|&\le C \delta \| \te +q_2  \|\,  \|\nabla v\| \le C \delta \left(\|\te\|^2+1\right) +  \frac{\delta\nu_1}{4} \|\nabla v\|^2,\\
	|I_7|&\le \delta \|h_1\|\, \| v\|\le C \delta   \| v\|.
\end{align*} 
   Integrating by parts   and using the condition   $\int_\T\diver v\,\dd z =0$, we get
\begin{align*}
  I_5 &=-\delta \int_{\T^2} \nabla p_s(x,y) \left(\int_\T   v(x,y,z)\, \dd z\right) \dd x\, \dd y  \\&= \delta\int_{\T^2}   p_s(x,y)   \left(\int_\T \diver v(x,y,z) \,\dd z\right) \dd x\, \dd y=0.
\end{align*}
To estimate $I_2$ and $I_3$, we   note that 
\begin{align*}
\left\lag\int_0^z \diver  (v+q_1) \,\dd \z \,
\p_z v,  v \right\rag&=\frac12 \int_{\T^3} \int_0^z \diver  (v+q_1) \,\dd \z \,
\p_z |v|^2 \dd x\, \dd y\,\dd z\\&= - \frac12 \int_{\T^3}   \diver  (v+q_1)   \,
  |v|^2 \dd x\, \dd y\,\dd z\\&= \left\lag \lag (v+q_1) , \nabla\rag v ,  v \right \rag . 
\end{align*}
Thus
\begin{align*}
 |I_2+I_3|&\le \delta\left(\|v\|+ \|q_1\| \right) \|\nabla q_1\|_\ty \|v\| + C \delta\left(\|\nabla v\|+\|\nabla q_1\| \right) \|\p_z q_1\|_\ty \|v\|\\&\le C\delta\left( \| v\|^2+1\right)+ \frac{\delta\nu_1}{4} \|\nabla v\|^2 .
 \end{align*}
Combining the estimates for $I_i$ with inequality \eqref{6.6}, we obtain
\begin{align}\label{6.7}
  \frac{\dd }{\dd t} \|  v\|^2+&\delta \nu_1\|\nabla  v\|^2+\delta\mu_1\|\p_{z} v\|^2 \le C\delta\left(\|  v\|^2+\|  \te\|^2+1\right). 	
\end{align}

 {\it Step~2.~$L^2$-estimate for $\te$.}  Now  we   take the scalar product in $L^2$ of Eq.~\eqref{6.5} with~$ \te $:    
    \begin{align}
     \frac{1}{2}\frac{\dd }{\dd t}  \|\te\|^2+& \delta\nu_2 \|\nabla  \te\|^2+\delta\mu_2\|\p_{z}\te\|^2= 
      -\delta \left\lag L_2  q_2 ,   \te\right\rag 
   -\delta  \left\lag \lag v +q_1 ,  \nabla\rag (\te + q_2),  \te\right\rag  \nonumber \\   &  +   \delta  \left\lag     \int_0^z \diver (v+q_1) \,\dd \z\,
\p_z (\te+q_2)   , \te\right\rag -  \left \lag \lag v ,  \nabla\rag \zeta ,  \te\right\rag   
   \nonumber   \\ &+   \left\lag  \int_0^z \diver  v  \,\dd \z\,
\p_z   \zeta    ,  \te\right\rag      +    \delta \left\lag  h_2,  \te\right\rag= \sum_{i=1}^6 J_i.  \label{6.8}
    \end{align}  We start with the terms $J_1, J_4, J_5, J_6$: 
      \begin{align*}
  	|J_1|&\le    \delta  \|L_2 q_2\|\,  \| \te\|\le    C\delta    \| \te\|, \\
  	|J_4|&\le \|v\| \, \|\nabla \zeta\|_\ty \|\te\|\le C  \|v\|  \, \|\te\|, \\
  	|J_5|&\le C \|\nabla v\|\, \|\p_z \zeta\|_\ty  \|\te \|\le C \|\nabla v\|\, \|\te \|,  \\
  		|J_6|&\le   \delta \|h_2\|\,\|  \te\| \le  C\delta \|   \te\| . 
  \end{align*}
To estimate $J_2$ and $J_3$, we use the equality 
\begin{align*}
\left\lag\int_0^z \diver  (v+q_1) \,\dd \z \,
\p_z \te,  \te \right\rag&=\frac12 \int_{\T^3} \int_0^z \diver  (v+q_1) \,\dd \z \,
\p_z (\te^2)\,\dd x\, \dd y\,\dd z\\&= - \frac12 \int_{\T^3}   \diver  (v+q_1)   \,
   (\te^2)\,\dd x\, \dd y\,\dd z\\&= \left\lag \lag (v+q_1) , \nabla\rag \te  ,  \te  \right \rag . 
\end{align*}
Then 
\begin{align*}
 |J_2+J_3|&\le \delta\left(\|v\|+ \|q_1\| \right)\|\nabla q_2\|_\ty \|\te\| + C \delta \left(\| \nabla v\|+\|\nabla q_1\|\right)  \|\p_z q_2\|_\ty  \|\te\| \nonumber\\&\le C\delta\left( \| v\|^2+\| \te\|^2+1\right)+  \delta\nu_1   \|\nabla v\|^2 .
 \end{align*}
The estimates for  $J_i$ and  \eqref{6.8} imply that
\begin{align*}
	 \frac{\dd }{\dd t}  \|\te\|^2+  \delta\nu_2 \|\nabla  \te\|^2+\delta\mu_2\|\p_{z}\te\|^2&\le C \left( \| v\|^2+\|\nabla  v\|^2+\| \te\|^2+1\right)\\&\quad + C\delta\left( \| v\|^2+\| \te\|^2+1\right) + \delta\nu_1  \|\nabla v\|^2 .
\end{align*}
Combining this with \eqref{6.7}, we get
\begin{gather}
		 \frac{\dd }{\dd t}  \|\te\|^2+\left(C(\nu_1\de )^{-1}+1\right) \frac{\dd }{\dd t}  \|v\|^2+  \delta\nu_2 \|\nabla  \te\|^2+\delta\mu_2\|\p_{z}\te\|^2\nonumber\\ \le C(1+\de) \left( \| v\|^2 +\| \te\|^2+1\right).\label{6.9}
\end{gather}
Integrating in time, we obtain 
 $$
  \|\te(t)\|^2+ \|v(t)\|^2 \le  C(1+\de)\int_0^t  \left( \| v\|^2 +\| \te\|^2+1\right)\dd s.
 $$The Gronwall inequality implies  
 $$
 \sup_{t\in [0,1],\, \de\in (0,1]} \left( \|\te(t)\|^2+ \|v(t)\|^2\right)<  \ty. 
 $$Going back to \eqref{6.7} and  \eqref{6.9}, we see that       
 \begin{align}\label{6.10}
 \|  v(t)\|^2+\|  \te(t)\|^2+ \int_0^t\left(   \|  v\|_1^2 +  \|  \te\|_1^2\right)\dd s\le    C\delta \quad \text{for   $t\in [0,1]$}. 
\end{align}

  {\it Step~3.~Barotropic-baroclinic formulation.} Following~\cite{CT-2007}, next we   
  estimate the~$L^6$-norm of $v$ and $H^1$-norms of $v$ and $\te$ by using barotropic-baroclinic formulation of PEs. More precisely, we denote
  $$ 
  \bar \phi =(2\pi)^{-1}\int_\T \phi(x,y,z)\,\dd z, \quad \tilde \phi=\phi-\bar \phi.
  $$Then the  barotropic mode $\bar v$ satisfies the following system of equations:
\begin{gather}
 \p_t \bar v   - \delta \nu_1 \Delta (\bar v+ \bar q_1)+ \delta\lag \bar v+\bar q_1, \nabla\rag (\bar v+\bar q_1) \nonumber\\+\delta \overline{\lag \tilde v+\tilde q_1, \nabla\rag (\tilde v+\tilde q_1)+ \diver (\tilde v+\tilde q_1) (\tilde v +\tilde q_1) }   +   \delta f   (\bar v+ \bar q_1)^\bot  \nonumber \\ +
 \delta \nabla p_s -\delta \overline{\int_0^z \nabla(\te +q_2)\,\dd \z}  = \de \bar h_1, \quad    
 \diver \bar v=0, \label{6.11}  
 \end{gather}and the baroclinic mode  $\tilde v$ the following one:
\begin{gather}
  \p_t \tilde v   + \delta L_1   (\tilde v+ \tilde q_1)+ \delta\lag \tilde v+\tilde q_1, \nabla\rag (\tilde v+\tilde q_1)-\delta\int_0^z \diver  (\tilde v +\tilde q_1)\,\dd \z \,
\p_z (\tilde v+\tilde q_1) 
   \nonumber\\+\delta\lag \bar  v+\bar q_1, \nabla\rag (\tilde v+\tilde q_1)+\delta\lag \tilde v+\tilde q_1, \nabla\rag (\bar v+\bar q_1) 
   \nonumber\\-\delta \overline{\lag \tilde v+\tilde q_1, \nabla\rag (\tilde v+\tilde q_1)+\diver (\tilde v+\tilde q_1) (\tilde v +\tilde q_1) } 
    +   \delta f   (\tilde v+ \tilde q_1)^\bot \nonumber\\-\delta {\int_0^z \nabla(\te +q_2)\,\dd \z}      +\delta \overline{\int_0^z \nabla(\te +q_2)\,\dd \z}  = \de \tilde h_1 \label{6.12}  
 \end{gather}(see \cite{CT-2007} for details). The advantage of this representation is that there is no pressure term in Eq.~\eqref{6.12} and the barotropic mode depends only on horizontal variables $(x,y)$ (its properties are similar to the ones   of 2D NS system).
 
  {\it Step~4.~$L^6$-estimate for $\tilde v$.}  We  take the scalar product in $L^2$ of Eq.~\eqref{6.12} with~$\tilde v |\tilde v|^4$:
 \begin{align}\label{6.13}
 \frac{1}{6}  \frac{\dd }{\dd t} \| \tilde v\|_{L^6}^6& +\delta \nu_1\||\nabla  \tilde v|\, |\tilde v|^2\|^2+\delta \nu_1\|\tilde v\,|\nabla  |\tilde v|^2|  \|^2\nonumber\\&+\delta \mu_1\||\p_z  \tilde v|\, |\tilde v|^2\|^2+\delta \mu_1\|\tilde v \,|\p_z  |\tilde v|^2|  \|^2 =  - \delta\left\lag L_1  \tilde q_1 , \tilde v |\tilde v|^4\right\rag\nonumber\\& 
 - \delta\left\lag\lag \tilde v+\tilde q_1, \nabla\rag (\tilde v+\tilde q_1)-\int_0^z \diver  (\tilde v +\tilde q_1)\,\dd \z \,
\p_z (\tilde v+\tilde q_1)  ,  \tilde v |\tilde v|^4\right\rag  \nonumber\\
&-\delta\left\lag \lag \bar  v+\bar q_1, \nabla\rag (\tilde v+\tilde q_1), \tilde v |\tilde v|^4 \right \rag 
-\delta\left\lag  \lag \tilde v+\tilde q_1, \nabla\rag (\bar v+\bar q_1), \tilde v |\tilde v|^4 \right \rag
\nonumber\\&+
\delta\left\lag\overline{\lag \tilde v+\tilde q_1, \nabla\rag (\tilde v+\tilde q_1)+\diver (\tilde v+\tilde q_1) (\tilde v +\tilde q_1) }   , \tilde v |\tilde v|^4 \right \rag  \nonumber\\&+ 
\delta\left\lag  {\int_0^z \nabla(\te +q_2)\,\dd \z}      - \overline{\int_0^z \nabla(\te +q_2)\,\dd \z}, \tilde v |\tilde v|^4 \right \rag\nonumber\\
&- 
\delta\left\lag f   (\tilde v+ \tilde q_1)^\bot, \tilde v |\tilde v|^4 \right \rag  + \delta \left\lag h_1 ,  \tilde v |\tilde v|^4 \right\rag = \sum_{i=1}^{8}I_i.
\end{align} Then
\begin{align*}
|I_1|&\le \delta\|L_1 \tilde q_1\|_{L^6}\|\tilde v\|_{L^6}^5	\le  C \delta \|\tilde v\|_{L^6}^5,	
\\|I_7|&= \delta |\lag f   \tilde q_1^\bot,  \tilde v |\tilde v|^4 \rag|\le  C\delta \|\tilde q_1\|_{L^6}  \|\tilde v\|_{L^6}^5	\le  C \delta \|\tilde v\|_{L^6}^5,
\\|I_8|&\le \delta\|\tilde h_1\|_{L^6}\|\tilde v\|_{L^6}^5	\le  C \delta \|\tilde v\|_{L^6}^5.
\end{align*}
Integrating by parts, we see that  
  $$
     \left\lag \lag \tilde v+\tilde q_1, \nabla\rag  \tilde v -\int_0^z \diver  (\tilde v +\tilde q_1)\,\dd \z \,
\p_z  \tilde v    , \tilde v |\tilde v|^4\right\rag  =0, 
  $$which implies   (again by integrating by parts)
  \begin{align*}
  	|I_2|&\le \delta \left(\|\tilde v\|_{L^6} +\|\tilde q_1\|_{L^6} \right) \|\nabla  \tilde q_1\|_\ty \|\tilde v\|_{L^6}^5 \\&\quad+C\delta \left(\|\tilde v\|_{L^6} +\|\tilde q_1\|_{L^6} \right) \|\nabla  \tilde q_1\|_\ty \|\tilde v\|_{L^6}^5 \\&\quad
  	 +C\delta  \left(   \|\tilde v\|_{L^6}
+ \| \tilde q_1\|_{L^6}  \right)  \|\tilde v\|_{L^6}^2 \|\p_z  \tilde q_1\|_\ty \||\nabla  \tilde v|\, |\tilde v|^2\|
\\&\le   C \delta    \left(\|\tilde v\|_{L^6}^6+1\right) +\frac{\delta\nu_1}{9}\||\nabla  \tilde v|\, |\tilde v|^2\|^2 .
  \end{align*}
As $ \diver \bar v= \diver \bar q_1=0$, we have that
  $$
     \left  \lag \lag \bar  v+\bar q_1, \nabla\rag  \tilde v   , \tilde v |\tilde v|^4\right\rag  =0, 
  $$
 hence 
 $$
 |I_3|\le \delta \left(\|\bar v\|_{L^6}+\|\bar q_1\|_{L^6} \right) \|\nabla \tilde q_1\|_\ty \|\tilde v \|_{L^6}^5 \le C\delta \left(\|\bar v\|_{L^6}+1\right)   \|\tilde v \|_{L^6}^5.
 $$
 To estimate $I_4$, we first integrate by parts:
 $$
 I_4=   \delta   \left\lag ( \bar v+\bar q_1)
 \diver \left( \tilde v+\tilde q_1\right), \tilde v|\tilde v|^4 \right\rag + \delta  \left\lag \lag  \tilde v+\tilde q_1, \nabla\rag (\tilde v |\tilde v|^4),  \bar v+\bar q_1\right\rag. 
 $$ We decompose $I_4$ as  $ I_4=I_4^1+I_4^2$, where  
 $$
 I_4^1=-\delta \left\lag \lag \tilde v ,\nabla\rag  \bar v  , \tilde v|\tilde v|^4\right\rag= 
  \delta   \left\lag   \bar v
 \diver \tilde v , \tilde v|\tilde v|^4 \right\rag + \delta  \left\lag \lag  \tilde v, \nabla\rag (\tilde v |\tilde v|^4),  \bar v \right\rag.
 $$It is proved on pages 255--257 in \cite{CT-2007} that
\begin{align*}
  |I_4^1|&\le C\delta \left(\|\bar v\|^{\frac12}\|\nabla\bar v\|^{\frac12}\|\tilde  v\|_{L^6}^{\frac32}\||\nabla  \tilde v|\, |\tilde v|^2\|^{\frac32}+
  \|\bar v\|^{\frac12}\|\nabla\bar v\|^{\frac12}\|\tilde  v\|_{L^6}^6
\right)
  \\&\le C\delta \left(\|\bar v\|^2\|\nabla\bar v\|^2+1\right)\|\tilde  v\|_{L^6}^6 +    \frac{\delta \nu_1}{9}  \||\nabla  \tilde v|\, |\tilde v|^2\|^2.
\end{align*}Then, integrating by parts,  we  estimate the   term $ I_4^2$ as follows:  
  \begin{align*}
 |I_4^2|&\le C \de\big(\|\bar v\|_{L^6} \|\nabla \tilde q_1\|_\ty\|\tilde v\|_{L^6}^5+ \|\bar q_1\|_{\ty}  \||\nabla  \tilde v|\, |\tilde v|^2\|\, \|\tilde v\|_{L^6}^3\\&\quad+  \|\bar q_1\|_\ty\|\nabla \tilde q_1\|_{L^6} \|\tilde v\|_{L^6}^5  +         \|\tilde q_1\|_{\ty} \|\tilde v\|_{L^6}^2  \||\nabla  \tilde v|\, |\tilde v|^2\|   \left( \|\bar v\|_{L^6}+\|\bar q_1\|_{L^6}\right)\big)\\&\le  C\de \left(\|\bar v\|_{L^6}  \|\tilde v\|_{L^6}^5+\|\bar  v\|_{L^6}^2\|\tilde  v\|_{L^6}^4+\|\tilde  v\|_{L^6}^6 +1\right)  +
 \frac{\delta \nu_1}{9}  \||\nabla  \tilde v|\, |\tilde v|^2\|^2   .  
  \end{align*}     To estimate $I_5$, we first integrate by parts:
\begin{align*}
I_5&=\delta \left\lag \overline{\lag \tilde v+\tilde q_1, \nabla\rag (\tilde v+\tilde q_1) 
   +   \diver (\tilde v+\tilde q_1) (\tilde v +\tilde q_1) }, \tilde v|\tilde v|^4 \right\rag \\&=-\delta\sum_{k,j=1}^2\int_{\T^3} \overline {(\tilde v+\tilde g_1)_k (\tilde v+\tilde g_1)_j} \p_k (\tilde v_j |\tilde v|^4),
\end{align*}where $\p_1=\p_x$ and $\p_2=\p_y$. We write $I_5=I_5^1+I_5^2$, where
$$
I_5^1=-\delta \sum_{k,j=1}^2\int_{\T^3} \overline {\tilde v_k \tilde v_j} \p_k (\tilde v_j |\tilde v|^4).
$$By the computations on pages 255--257 in \cite{CT-2007},    we have
\begin{align*}
|I_5^1|&\le C\delta \|\tilde  v\|_{L^6}^3\left(\|\tilde v\|+ \|\nabla\tilde v\|\right)\, \||\nabla  \tilde v|\, |\tilde v|^2\|\\&\le C\delta \|\tilde  v\|_{L^6}^6 \left(\|\tilde v\|^2+ \|\nabla\tilde v\|^2\right) +  \frac{\delta \nu_1}{9}  \||\nabla  \tilde v|\, |\tilde v|^2\|^2. \end{align*}
 Then we estimate $I_5^2$ as follows: 
\begin{align*}
|I^2_5|&\le C\delta\int_{\T^2}\left( \left(\int_\T \left(|\tilde v|+1\right) \,\dd z\right)^2 \int_\T |\nabla \tilde v||\tilde v|^4 \,\dd z\right) \dd x\,\dd y\\&\le 
  C\delta\int_{\T^2}\left(\int_\T |\tilde v|^2 \,\dd z \int_\T |\nabla \tilde v||\tilde v|^4 \,\dd z\right) \dd x\,\dd y+	C\delta \int_{\T^3} |\nabla \tilde v||\tilde v|^4 \dd x\,\dd y\,\dd z. 
\end{align*}The first term on the right-hand side is estimated exactly in the same way as $I_5^1$ (see \cite{CT-2007}), and the second term by
\begin{align*}
C\delta \int_{\T^3} |\nabla \tilde v||\tilde v|^4 \dd x\,\dd y\,\dd z\le C\delta  \||\nabla  \tilde v|\, |\tilde v|^2\| \,\|\tilde v\|_{L^4}^2\le C\de \|\tilde v\|_{L^4}^4+  \frac{\delta \nu_1}{9}  \||\nabla  \tilde v|\, |\tilde v|^2\|^2. 
\end{align*}
It remains to estimate $I_6$. To this end, we write $I_6=I_6^1+I_6^2$, where 
  $$
  I_6^1=\delta \left\lag    {\int_0^z \nabla \te  \,\dd \z}      - \overline{\int_0^z \nabla \te  \,\dd \z}, \tilde v |\tilde v|^4 \right\rag.
  $$
Again we refer to  pages 255--257 in \cite{CT-2007} for the proof of the following inequality:
 \begin{align*}
 |I_6^1|&\le 	C\delta \|\bar \te\|^{\frac12} \|\nabla\bar \te\|^{\frac12} \|\tilde v\|_{L^6}^\frac32  \left(\|\tilde v\|^{\frac12} +\|\nabla\tilde v\|^{\frac12}\right)   \||\nabla  \tilde v|\, |\tilde v|^2\| \\&\le C\delta \left(\|\bar \te\|^2  \|\nabla\bar \te\|^2+ \|\tilde v\|_{L^6}^6 \|\tilde v\|^2+  \|\tilde v\|_{L^6}^6 \|\nabla\tilde v\|^2\right)   + \frac{\delta \nu_1}{9}  \||\nabla  \tilde v|\, |\tilde v|^2\|^2.
  \end{align*}
 Then we estimate $I_6^2$ as follows:
   \begin{align*}
 |I_6^2|&\le C\delta \|\nabla q_2\|_\ty \|\tilde v\|_{L^5}^5\le C\delta   \|\tilde v\|_{L^5}^5.
 \end{align*}
   Combining all the above estimates for the terms $I_i$ with \eqref{6.13}, \eqref{6.10}, the Sobolev embedding $H^1(\T^3)\subset L^6(\T^3)$, and the Gronwall inequality, we conclude~that 
   \begin{equation}\label{6.14}
   \|\tilde v(t)\|_{L^6}+\int_0^t\||\nabla  \tilde v|\, |\tilde v|^2\|^2\,\dd s\le C\delta	\quad \text{for   $t\in [0,1]$.}
   \end{equation} 

  {\it Step~5.~Estimate for $\nabla \bar v$.} Here we take the scalar product in $L^2$ of the first equation in~\eqref{6.11} with $-\Delta \bar v$:
    \begin{align}\label{6.15}
\frac{1}{2}  \frac{\dd }{\dd t} \|\nabla \bar v\|^2+&\delta \nu_1\|\Delta \bar v\|^2=   -\delta \nu_1 \lag \Delta \bar q_1 , \Delta \bar v\rag 
+\delta \left \lag  \lag \bar v+\bar q_1, \nabla\rag (\bar v+\bar q_1)  , \Delta \bar v \right\rag
\nonumber\\&
+\delta \left \lag  \overline{\lag \tilde v+\tilde q_1, \nabla\rag (\tilde v+\tilde q_1)+ \diver (\tilde v+\tilde q_1) (\tilde v +\tilde q_1) }  , \Delta \bar v \right\rag
\nonumber\\&+\delta \left \lag f   (\bar v+ \bar q_1)^\bot   , \Delta \bar v \right\rag
 +\delta \left \lag  \nabla p_s  , \Delta \bar v \right\rag
\nonumber\\&-\delta \left \lag \overline{\int_0^z \nabla(\te +q_2)\,\dd \z}   , \Delta \bar v \right\rag - \delta \left\lag \bar h_1 , \Delta \bar v \right\rag = \sum_{i=1}^7 I_i.
\end{align}We estimate\footnote{In the estimate for $I_4$, we use the equality $\lag \bar v^\bot, \Delta \bar  v\rag=0$ which is easily verified by integration by parts in $x$ and $y$.} $I_1, I_4, I_7$
as follows:
\begin{align*}
|I_1|&\le \delta \nu_1\| \Delta \bar q_1\|\, \|\Delta\bar  v\| 	\le  C \delta + \frac{\delta \nu_1}{10} \| \Delta \bar  v\|^2,	
\\|I_4|&= \delta |\lag f   \bar q_1^\bot,  \Delta \bar  v   \rag|\le C \delta + \frac{\delta \nu_1}{10} \| \Delta \bar  v\|^2,
\\|I_7|&\le \delta\|\bar h_1\|\, \|\Delta \bar  v\| 	\le  C \delta + \frac{\delta \nu_1}{10} \| \Delta \bar  v\|^2.
\end{align*}
Integrating by parts  and using the fact that $\diver \bar v=0$, we get~$  I_5=I_6   =0.$
Next we use the H\"older inequality, the  Sobolev embedding~$H^1(\T^2)\subset L^4(\T^4)$, and  the interpolation inequality to estimate   $I_2$:
\begin{align*}
|I_2|&\le \delta \| \bar v+\bar q_1\|_{L^4} \| \nabla(\bar v+\bar q_1)\|_{L^4}	 \|\Delta \bar v\|\le \delta \| \bar v+\bar q_1\|^\frac12 \| \nabla(\bar v+\bar q_1)\|\, 	 \|\Delta \bar v\|^\frac32\\
&\le C\delta   \left(\| \nabla \bar v\|^6+1\right) + \frac{\delta \nu_1}{10} 	 \|\Delta \bar v\|^2.
\end{align*}
Finally, we use the H\"older inequality to estimate $I_3$:
\begin{align*}
|I_3|&\le C\delta \int_{\T^2}\left(\int_\T |\tilde v+\tilde q_1|\,|\nabla (\tilde v+\tilde q_1)|\,\dd z\right) |\Delta \bar v|\,\dd x\,\dd y
\\&\le C\delta \int_{\T^2}\left(\int_\T (|\tilde v|+1)\,(|\nabla \tilde v|+1)\,\dd z\right) |\Delta \bar v|\,\dd x\,\dd y
\\&\le C\delta \int_{\T^2}\left(\int_\T (|\nabla  \tilde v|+ 1) \,\dd z\right)^\frac12 \left(\int_\T (|\tilde v|+ 1)^2(|\nabla  \tilde v|+1)\,\dd z\right)^\frac12 |\Delta \bar v|\,\dd x\,\dd y\\&\le  C\delta \left(\|\nabla \tilde v\|^2+\|\tilde v\|_{L^4}^2  + \||\nabla  \tilde v|\, |\tilde v|^2\|^2 +1\right)       +\frac{\delta \nu_1}{10} 	 \|\Delta \bar v\|^2.	
\end{align*}
 Combining the estimates for $I_i$ and inequalities \eqref{6.15}, \eqref{6.14}, and \eqref{6.10}, we~get
\begin{equation}\label{6.16}
   \frac{\dd }{\dd t} \|\nabla \bar v\|^2+ \delta \nu_1\|\Delta \bar v\|^2\le C\de \left(\| \nabla \bar v\|^6+1\right) .
\end{equation}
This inequality implies that 
\begin{equation}
	\label{6.17}
     \|\nabla \bar v(t)\|^2+ \int_0^t  \|\Delta \bar v\|^2 \dd s\le C\de  \quad \text{for }t\in [0,1],   
\end{equation} provided that   $\delta>0$ is sufficiently small.
Indeed, to see this, let   $A_\delta=C \delta  $  and
$$
		\Phi(t)= A_\delta + A_\delta\int_0^t  \|\nabla \bar v\|^6 \,\dd s.   
		  $$
		  By inequality~\eqref{6.16}, we have   
$$
		\left(\frac{\dd }{\dd t} \Phi\right)^{\frac13}\le A_\delta^{\frac13}\Phi,
$$
	or, equivalently,
	$$
		\frac{1}{\Phi^{3}} \frac{\dd }{\dd t} \Phi \le A_\delta.
	$$ 
	Integrating this, we obtain
	$$
		\Phi(t)\le  A_\delta \left(1-2A_\delta^{3} t\right)^{-\frac12} \quad \text{for  $0 \le t <1  \wedge \left( \frac1{2A_\delta^{3}}\right) $}.
	$$
	Thus 
	 $$
		\Phi(t)\le 2  A_\delta   \quad \text{for  $0 \le t <1  \wedge \left(\frac{3}{8A_\delta^{3}}\right)$}.
	$$  	Choosing $\delta_0>0$ so small that 
$$
\frac{3}{8A_\delta^{3}}>1 \quad\text{for $ \de\in (0,\de_0)$}, 
$$we arrive at
$$
\Phi(t)\le 2  A_\delta\quad  \text{for  $  t \in [0,  1]$, $   \delta\in (0,\de_0)$}. 
$$
Combining this with~\eqref{6.16}, we prove  \eqref{6.17}. Below everywhere we shall assume that $\de\in (0,\de_0)$.

  {\it Step~6.~Estimate for $\p_zv$.} 
  The function $\omega=\p_z v$ is a solution of the equation
  \begin{gather}
  \p_t\omega  + \delta L_1(\omega+\hat q_1)+ \delta\lag v+q_1, \nabla\rag (\omega+\hat q_1)    - \delta\int_0^z \diver  (v +q_1)\,\dd \z \,
\p_z (\omega+\hat q_1)  \nonumber\\
+ \delta\lag \omega+\hat q_1, \nabla\rag (v+q_1)    - \delta  \diver  (v +q_1)  \,
  (\omega+ \hat q_1)  \nonumber\\+   \delta f   (\omega+\hat q_1)^\bot    -\delta   \nabla(\te +q_2)   = \de \hat h_1,  \label{6.18} 
 \end{gather}where we denote $\hat q_1=\p_z q_1$ and $\hat h_1=\p_z h_1$.
 Let us take the scalar product in $L^2$ of Eq.~\eqref{6.18} with~$\omega$:
  \begin{align}\label{6.19}
 \frac{1}{2} & \frac{\dd }{\dd t} \|\omega\|^2+\delta \nu_1\|\nabla \omega\|^2+\delta\mu_1\|\p_{z} \omega\|^2 =   -\delta\left\lag L_1  \hat q_1 , \omega\right\rag \nonumber\\ &-\delta\left\lag \lag v+q_1, \nabla\rag (\omega+\hat q_1)    - \int_0^z \diver  (v +q_1)\,\dd \z \,
\p_z (\omega+\hat q_1) , \omega \right \rag   \nonumber\\
&- \delta\left\lag \lag \omega+\hat q_1, \nabla\rag (v+q_1)    -   \diver  (v +q_1)  \,
  (\omega+ \hat q_1), \omega \right\rag   \nonumber   \\&  -  \delta \left\lag f   (\omega+\hat q_1)^\bot, \omega\right\rag      +\delta\left\lag  \nabla(\te +q_2)  , \omega \right\rag + \delta \lag \hat h_1 , \omega \rag = \sum_{i=1}^6I_i.
\end{align} Then $I_1, I_4, I_5, I_6$ are estimated as follows:
 \begin{align*}
|I_1|&\le \delta \|L_1 q_1 \| \, \| \omega\|\le  C \delta \|\omega\|,   \\
	| I_4|&= \delta |\lag f   \hat q_1^\bot, \omega \rag|\le     \delta \|f\hat q_1\|\, \|\omega\|\le   C\delta \|\omega\|,      \\
	|I_5|&\le C \delta \| \te +q_2  \|\,  \|\nabla \omega\| \le C \delta \left(\|\te\|^2+1\right) +  \frac{\delta\nu_1}{4} \|\nabla \omega\|^2,\\
	|I_6|&\le \delta \|\hat h_1\|\, \|\omega\|\le C \delta \|\omega\| .
\end{align*} 
Integrating by parts in $z$,   we get
$$
\left\lag \lag v+q_1, \nabla\rag  \omega     - \int_0^z \diver  (v +q_1)\,\dd \z \,
\p_z  \omega  , \omega \right \rag=0,
$$so we can estimate $I_2$ by
\begin{align*}
|I_2|&\le 	\delta \left(\|v\|+\|q_1\|\right) \|\nabla \hat q_1\|_\ty \|\omega\| +C\delta \left(\|\nabla v\|+\|\nabla q_1\|\right) \| \p_z \hat q_1\|_\ty \|\omega\|
\\&\le C\delta \left(\|v\|_1+1\right) \|\omega\|.
\end{align*}Integrating by parts in $x$ and $y$  and
 using the H\"older,  Gagliardo--Nirenberg, and Sobolev inequalities,   we obtain
\begin{align*}
|I_3|&\le 	C\delta  \int_{\T^3}  |v+q_1| \left( |\nabla (\omega+\hat q_1)| |\omega|+|  \omega+\hat q_1 | |\nabla \omega|\right)\dd  x\, \dd y\, \dd z \\&\le 	C\delta  \int_{\T^3}  \left(|v|+1\right)   \left(|\omega|+1\right) \left( |\nabla  \omega|+1\right) \dd  x\, \dd y\, \dd z\\&\le 	C\delta  \left(\|v\|_{L^6}+1\right)\left(\|\omega\|_{L^3}+1\right)\left(\|\nabla \omega\|+1\right)
\\&\le 	C\delta  \left(\|v\|_{L^6}+1\right)\left(\|\omega\|^\frac12\|\omega\|_1^\frac12+1\right)\left(\|\nabla \omega\|+1\right)
\\&\le 	C\delta  \left( \|v\|_{L^6}^4\|\omega\|^2+\|v\|_{L^6}^2+\|\omega\|^2+1 \right)+\frac{\delta \nu_1}{2}\|\nabla \omega\|^2+\frac{\delta\mu_1}{2}\|\p_{z} \omega\|^2
\\&\le 	C\delta  \left( \left(\|\nabla \bar v\|^4+\|\tilde  v\|_{L^6}^4\right)\|\omega\|^2+\|\nabla \bar v\|^2+\| \tilde v\|_{L^6}^2+\|\omega\|^2+1 \right)\\&\quad +\frac{\delta \nu_1}{4}\|\nabla \omega\|^2+\frac{\delta\mu_1}{2}\|\p_{z} \omega\|^2.
\end{align*} Combining the estimates for $I_i$ with  \eqref{6.19}, the Gronwall inequality, and inequalities \eqref{6.10},~\eqref{6.14}, and~\eqref{6.17}, we obtain 
\begin{equation}\label{6.20}
     \|\p_z v (t)\|^2+  \int_0^t\| \p_zv\|_1^2\dd s   \le C\delta  \quad \text{for $t\in[0,1]$}  .
     \end{equation}

   {\it Step~7.~Estimate for $\nabla v$.}~We   take the scalar product in $L^2$ of Eq.~\eqref{6.4} with~$-\Delta v$:
  \begin{align}\label{6.21}
\frac{1}{2}  \frac{\dd }{\dd t} \|\nabla v\|^2+&\delta \nu_1\|\Delta v\|^2+\delta\mu_1\|\nabla\p_{z} v\|^2=   \delta\lag L_1 q_1 , \Delta v\rag \nonumber\\ &+\delta\left\lag \lag v+q_1 , \nabla\rag (v+q_1 ), \Delta v \right \rag   \nonumber\\
&- \delta\left\lag\int_0^z \diver  (v +q_1 )\,\dd \z \,
\p_z (v+q_1 ), \Delta v \right\rag   \nonumber   \\&  +  \delta \left\lag f   (v+q_1 )^\bot, \Delta v\right\rag    +
 \delta\left\lag\nabla p_s, \Delta v\right\rag  \nonumber\\&-\delta\left\lag \int_0^z \nabla(\te +q_2)\,\dd \z , \Delta v \right\rag - \delta \left\lag h_1 , \Delta v \right\rag = \sum_{i=1}^7 I_i.
\end{align}
 Again the terms      $ I_1,  I_4,  I_6,  I_7$  are easier to estimate:\footnote{In the estimate for $ I_4$, we use the equality  $\lag v^\bot,  \Delta v\rag=0$.}
    \begin{align*}
| I_1|&\le \delta \|L_1 q_1 \| \, \|\Delta v\|\le  C \delta +      \frac{\delta\nu_1}{16} \|\Delta   v\|^2, \\
	|  I_4|&= \delta |\lag f    q_1^\bot, \Delta v \rag|\le     \delta \|fq_1\|\,   \|\Delta v\|\le   C\delta    +  \frac{\delta\nu_1}{16} \|\Delta v\|^2,   \\
	| I_6|&\le C\delta \|\te +q_2 \|_1\,  \|\Delta v\| \le C \delta \left(\|\te\|_1^2+1\right) +  \frac{\delta\nu_1}{16} \|\Delta v\|^2,\\
	| I_7|&\le \delta \|h_1\|\, \|\Delta v\|\le C \delta +\frac{\delta\nu_1}{16} \|\Delta v\|^2.
\end{align*} 
   Integrating by parts in     $x$ and $y$ and using the condition   $\int_\T\diver v\,\dd z =0$, we get
\begin{align*}
  I_5 &=\delta \int_{\T^2} \nabla p_s(x,y) \left(\int_\T \Delta v(x,y,z)\, \dd z\right) \dd x\, \dd y \\&=\delta\int_{\T^2} \nabla p_s(x,y)\, \Delta \left(\int_\T v(x,y,z) \,\dd z\right) \dd x\, \dd y \\&=-\delta\int_{\T^2}   p_s(x,y)\, \Delta \left(\int_\T \diver v(x,y,z) \,\dd z\right) \dd x\, \dd y=0.
\end{align*}
  We decompose the terms $  I_2$ and $  I_3$ as follows: 
\begin{align*}
	  I_2&= P_1+  P_2, \,\,\quad   P_1=\delta\left\lag \lag v, \nabla\rag v, \Delta v \right \rag,\\
    I_3&=   Q_1+Q_2, \quad   Q_1= -\delta\left\lag\int_0^z \diver  v  \,\dd \z \,
\p_z  v , \Delta v \right\rag,
\end{align*}
and estimate the quadratic in $v$ terms   $  P_2$ and $  Q_2$ in the following way:
\begin{align*}
|  P_2|&\le  \delta  \left(\|v\|\, \|\nabla q_1 \|_\ty + \|\nabla v\| \,\|q_1 \|_\ty+\| q_1\| \,\|\nabla q_1 \|_\ty \right)  \|\Delta  v\|
  \\ &\le  C\delta \left(\|v\|_1^2+1\right)+\frac{\delta\nu_1}{16} \|\Delta  v\|^2,\\
|  Q_2|&\le 	C  \delta  \left(\|\p_z v\|\, \|\nabla q_1 \|_\ty + \|\nabla  v\| \,\|\p_z q_1 \|_\ty+\|\nabla  q_1\| \,\|\p_z q_1 \|_\ty  \right)  \|\Delta v \|
 \\& \le	 C\delta \left (\|v\|_1^2+1\right)+\frac{\delta\nu_1}{16} \|\Delta  v\|^2.
\end{align*}   For $  P_1$, we use the H\"older, Sobolev, and Gagliardo--Nirenberg inequalities: 
\begin{align*}
|  P_1| &\le C \delta \int_{\T^3}   |v||\nabla v| |\Delta v|\,\dd x\,\dd y\,\dd z  \le  C\delta  \| v\|_{L^6}\| \nabla v\|_{L^3} \|\Delta  v\|\\&  
 \le C  \delta  \| v\|_{L^6} \| \nabla v\|^\frac12  \| \nabla v\|_1^\frac12  \|\Delta  v\|
\\&  
 \le  C \delta  \left(\| \nabla \bar v\|^4+\|\tilde v\|^4_{L^6}\right)\|\nabla v\|^2 + C\de  \|\nabla \p_z  v\|^2+\frac{\delta\nu_1}{16} \|\Delta  v\|^2.
\end{align*}
Next, we use the following inequality which is  proved in    Proposition~2.2 in~\cite{CT-2003}:
\begin{equation}\label{6.22} 
	\left|\left\lag\int_0^z \diver  \phi  \,\dd \z \,
   \varphi,  \psi \right\rag\right|\le C \|\phi\|_1^{\frac12}\|L_1 \phi\|^{\frac12}  \|\varphi\|^{\frac12}\|\varphi\|_1^{\frac12} \|\psi\| 
\end{equation}for any $\phi\in H^2(\T^3,\R^2)$, $\varphi\in H^1(\T^3,\R)$, and $\psi\in L^2(\T^3,\R)$.  Applying \eqref{6.22}, we obtain
\begin{align*}
| Q_1| & \le C \delta \| v\|_1^{\frac12}  \|L_1 v\|^{\frac12}  \|\p_z v\|^{\frac12}   \|\p_z v\|_1^{\frac12} \|\Delta v\|\\&\le C \delta \left(\|\nabla v\|^2 +\| \p_z v\|^2+\| v\|^2 \right)  \|\p_z v\|^2 \|\nabla \p_zv\|^2 +C\delta  \| \p_{zz}  v\|^2    +\frac{\delta\nu_1}{16} \|\Delta  v\|^2.
\end{align*}
The estimates for $ P_1,  P_2,  Q_1,  Q_2$ and $I_i$, together with \eqref{6.21}, the Gronwall inequality, and inequalities~\eqref{6.10},~\eqref{6.14},~\eqref{6.17}, and~\eqref{6.20}
imply that
$$\|\nabla v(t)\|^2+ \int_0^t \|\Delta v\|^2\dd s\le C\delta \quad\text{for $t\in [0,1]$}.	
$$Combining this with \eqref{6.20}, we get
\begin{equation}\label{6.23}
\| v(t)\|_1^2+ \int_0^t \|v \|_2^2\,\dd s\le C\delta \quad\text{for $t\in [0,1]$},	
\end{equation}which implies, in particular that $\|v(1)\|_1\to 0$ as $\delta\to 0^+$.

   {\it Step~8.~Estimate for $\p_z\te $.} Now we turn to the $H^1$-estimates for $\te$.
  To estimate~$\p_z\te $, we take the scalar product in $L^2$ of Eq.~\eqref{6.5} with $-\p_{zz} \te$:
      \begin{align}
     \frac{1}{2}\frac{\dd }{\dd t}  \|\p_z\te\|^2+& \delta\nu_2 \|\nabla \p_z\te\|^2+\delta\mu_2\|\p_{zz}\te\|^2= 
      \delta \left\lag L_2  q_2 , \p_{zz} \te\right\rag 
     \nonumber \\   &+\delta  \left\lag \lag v +q_1 ,  \nabla\rag (\te + q_2), \p_{zz}\te\right\rag+    \left \lag \lag v ,  \nabla\rag \zeta , \p_{zz}\te\right\rag \nonumber   \\ &-\!   \delta  \left\lag     \int_0^z \diver (v+q_1) \,\dd \z\,
\p_z (\te+q_2)   , \p_{zz}\te\right\rag   
-   \left\lag  \int_0^z \diver  v  \,\dd \z\,
\p_z   \zeta    , \p_{zz}\te\right\rag  \nonumber  \\  &-    \delta \left\lag  h_2, \p_{zz}\te\right\rag= \sum_{i=1}^6 J_i.  \label{6.24}
    \end{align}  We start with the terms $J_1$ and $J_6$: 
      \begin{align*}
  	|J_1|&\le    \delta  \|L_2 q_2\|\,  \|\p_{zz} \te\|\le C\de +  \frac{\delta\mu_2}{12}\|\p_{zz} \te\|^2, \\
  		|J_6|&\le   \delta \|h_2\|\,\|\p_{zz} \te\| \le  C\delta+ \frac{\delta\mu_2}{12}\|\p_{zz} \te\|^2. 
  \end{align*}
 To estimate $J_3$ and $J_5$, we integrate by parts in $z$ and use   the  Cauchy--Schwarz  and Sobolev  inequalities: 
\begin{align*}
  	|J_3|&\le | \left \lag \lag \p_z v ,  \nabla\rag \zeta   , \p_{z}\te\right\rag|+| \left \lag \lag v ,  \nabla\rag \p_z\zeta   , \p_{z}\te\right\rag|\\&\le  \| \p_z v\|\, \|\nabla \zeta\|_\ty \|\p_z\te\|+\|v\|\,   \|\nabla \p_z\zeta\|_{\ty}\|\p_{z}\te\|\le C  \|v\|_1 \|\p_z\te\|,\\
  	|J_5|&\le  | \lag   \diver  v   \,
\p_z   \zeta    , \p_{z}\te\rag  |+ \left| \bigg\lag  \int_0^z \diver  v  \,\dd \z\,
\p_{zz}   \zeta    , \p_{z}\te\bigg\rag\right|\\& \le C  \|\nabla  v\|\,   \|\p_z   \zeta\|_\ty    \| \p_{z}\te\|+ C\|\nabla  v\|\, \|\p_{zz}\zeta\|_\ty \|\p_z\te\| \le C  \|v\|_1 \|\p_z\te\|.
  \end{align*}
 We write the terms $J_2$ and $J_4$ as follows:
 \begin{align*}
	J_2&=P_1+P_2, \,\,\quad P_1=\delta\left\lag \lag v,  \nabla\rag \te , \p_{zz}\te\right\rag,\\
	J_4&=  Q_1+Q_2, \quad Q_1= -\delta\left\lag\int_0^z \diver  v  \,\dd \z \,
\p_z  \te  , \p_{zz}\te  \right\rag,
\end{align*}and estimate $P_2$ and $Q_2$ as in the previous steps: 
 \begin{align*}
| P_2|&\le  \delta  \left(\|v\|\, \|\nabla q_2 \|_\ty +  \|\nabla \te\|\,\|q_1 \|_\ty +\| q_1\| \,\|\nabla q_2 \|_\ty \right)  \|\p_{zz}  \te\|
  \\ &\le  C\delta \left(\|v\|_1^2+\|\te\|_1^2+1\right)+\frac{\delta\mu_2}{12} \|\p_{zz}  \te\|^2,\\
| Q_2|&\le	C  \delta  \left(\|\nabla  v\|\, \| \p_z q_2 \|_\ty + \|\p_z \te\|\, \|\nabla   q_1\|_\ty  + \|\nabla   q_1\|\, \|\p_z q_2 \|_\ty   \right)  \|\p_{zz} \te \|
 \\& \le	 C\delta  \left(\|v\|_1^2+\|\te\|_1^2+1\right)+\frac{\delta\mu_2}{12} \|\p_{zz}  \te\|^2.
\end{align*}
 To estimate $P_1$, we use the H\"older, Sobolev, and Gagliardo--Nirenberg inequalities:
\begin{align*}
|P_1|&\le \de \|v\|_{L^6}  \|\nabla \te \|_{L^3}  \|\p_{zz}\te\| \le C \de \|v\|_1 \|\nabla \te \|^\frac12   \|\nabla \te \|^\frac12_1  \|\p_{zz}\te\| \\&\le    C\de \|v\|_1^4\| \te \|_1^2  + \delta\nu_2 \|\Delta  \te\|^2 +  \frac{\delta\mu_2}{12}\| \p_{zz}\te\|^2.
\end{align*}
To estimate $Q_1$, we first  integrate by parts in $z$: 
\begin{align*}
\left\lag\int_0^z \diver  v  \,\dd \z \,
\p_z  \te  , \p_{zz}\te  \right\rag&=\frac12\int_{\T^3} \int_0^z \diver  v  \,\dd \z\, \p_z((\p_z \te)^2)\,\dd  x\, \dd y\, \dd z\\&=-\frac12\int_{\T^3}  (\p_z \te)^2\,  \diver  v\,\dd  x\, \dd y\, \dd z,
\end{align*}
then we use  the Cauchy--Schwarz and Gagliardo--Nirenberg inequalities:
\begin{align*}
|Q_1|&\le C \delta  \| \nabla v\| \, \|\p_z\te\|_{L^4}^2 	\le C \delta \| \nabla v\|\, \|\p_z\te \|^{\frac12}  \|\p_z\te\|_1^{\frac32}\\&\le  C \delta \|  v\|_1^4  \|\te \|_1^2  +\frac{\delta\nu_2}{2}\|\nabla \p_z\te\|^2+\frac{\delta\mu_2}{12}\| \p_{zz}\te\|^2. \end{align*}
Combining the estimates for $P_1,P_2,Q_1,Q_2$ and $J_i$ with inequalities~\eqref{6.10},~\eqref{6.23}, and \eqref{6.24}, we obtain       
  \begin{equation}\label{6.25}
  	  \|\p_z\te(t)\|^2 +\delta\mu_2\int_0^t\|\p_{zz}\te\|^2\dd s\le C \de+ \delta\nu_2    \int_0^t\|\Delta  \te\|^2   \dd s.
  \end{equation}

         {\it Step~9.~Estimate for $\nabla \te $.} Finally, 
to estimate $\nabla\te $, we take the scalar product in $L^2$ of Eq.~\eqref{6.5} with $-\Delta \te$:
     \begin{align}
     \frac{1}{2}\frac{\dd }{\dd t} \|\nabla\te\|^2+&\delta\nu_2 \|\Delta \te\|^2+\delta\mu_2\|\nabla \p_{z}\te\|^2= 
      \delta \left\lag L_2  q_2 , \Delta \te\right\rag 
     \nonumber \\&  +\delta  \left\lag \lag v +q_1 ,  \nabla\rag (\te + q_2), \Delta\te\right\rag+    \left \lag \lag v ,  \nabla\rag \zeta , \Delta \te\right\rag \nonumber   \\ &-   \delta  \left\lag     \int_0^z \diver (v+q_1) \,\dd \z\,
\p_z (\te+q_2)   , \Delta\te\right\rag   
-   \left\lag  \int_0^z \diver  v  \,\dd \z\,
\p_z   \zeta    , \Delta\te\right\rag  \nonumber  \\  &-    \delta \left\lag  h_2, \Delta\te\right\rag= \sum_{i=1}^6  J_i. \label{6.26} 
    \end{align}
The terms $ J_1$ and $ J_6$ are estimated as follows:  
  \begin{align*}
  	| J_1|&\le   \delta  \|L_2 q_2\|\,  \|\Delta  \te\|\le C\de +  \frac{\delta\nu_2}{12}\| \Delta \te\|^2, \\
  		| J_6|&\le   \delta \|h_2\|\,\|\Delta \te\| \le  C\delta+ \frac{\delta\nu_2}{12}\|\Delta \te\|^2. 
  \end{align*}
    To estimate $ J_3$ and $ J_5$, we integrate by parts   and use   the  Cauchy--Schwarz and Sobolev  inequalities:  
\begin{align*}
  	| J_3|&\le     C  \|v\|_1  \|\te\|_1,\\
  	| J_5|&\le   C  \|\Delta v\|\, \| \te\|_1.
  \end{align*}
The terms $ J_2$  and $ J_3$ are decomposed as follows:
 \begin{align*}
	 J_2&= P_1+ P_2, \,\,\quad  P_1=\delta\left\lag \lag v,  \nabla\rag \te , \Delta\te\right\rag,\\
	 J_3&=   Q_1+ Q_2, \quad  Q_1= -\delta\left\lag\int_0^z \diver  v  \,\dd \z \,
\p_z  \te  , \Delta\te  \right\rag,
\end{align*}and $ P_2$ and $ Q_2$ are estimated by
 \begin{align*}
|  P_2|&\le  \delta  \left(\|v\|\, \|\nabla q_2 \|_\ty +  \|\nabla \te\|\,\|q_1 \|_\ty +\| q_1\| \,\|\nabla q_2 \|_\ty \right)  \|\Delta   \te\|
  \\ &\le  C\delta \left(\|v\|_1^2+\|\te\|_1^2+1\right)+\frac{\delta\nu_2}{12} \|\Delta   \te\|^2,\\
|  Q_2|&\le	C  \delta  \left(\|\nabla   v\|\, \| \p_z q_2 \|_\ty +  \|\p_z \te  \|\, \|\nabla   q_1\|_\ty+\|\nabla   q_1\| \,\|\p_z q_1 \|_\ty  \right)  \|\Delta  \te \|
 \\& \le	 C\delta  \left(\|v\|_1^2+\|\te\|_1^2+1\right)+\frac{\delta\nu_2}{12} \|\Delta   \te\|^2.
\end{align*}
 Using the H\"older, Sobolev, and Gagliardo--Nirenberg inequalities, we obtain
\begin{align*}
| P_1|&\le \de \|v\|_{L^6}  \|\nabla \te \|_{L^3}  \|\Delta \te\| \le C \de \|v\|_1  \|\nabla \te \|^\frac12   \|\nabla \te \|^\frac12_1  \|\Delta \te\| \\&\le  C \de \|v\|_1^4  \| \te \|_1^2 +\frac{\delta\mu_2}{4}\| \p_{zz}\te\|^2+  \frac{\delta\nu_2}{12}\| \Delta \te\|^2.
\end{align*}
By \eqref{6.22}, we have 
\begin{align*}
| Q_1| &\le C \delta \|v\|_1^{\frac12} \|L_1v\|^{\frac12} \|\p_z\te\|^{\frac12} \|\p_z\te\|_1^{\frac12}  \|\Delta \te\|\\&\le C\delta \|v\|_1^2 \|L_1v\|^2 \|\p_z\te\|^2 +\frac{\delta\mu_2}{4}\|\p_{zz}\te\|^2  + \frac{\delta\nu_2}{12}\|\Delta \te\|^2.	
\end{align*}
 The estimates of $ P_1, P_2, Q_1, Q_2$ and $I_i$
 and the inequalities~\eqref{6.10},~\eqref{6.23}, and~\eqref{6.26} imply that  
$$ \|\nabla\te(t)\|^2+\delta\nu_2 \int_0^t\|\Delta \te\|^2\,\dd s  \le  C\delta +\delta\mu_2\|\p_{zz}\te\|^2 .   
$$ From this and \eqref{6.25} we derive that $\|\te (t)\|_1\le C\delta$, so 
  $\|\te(1)\|_1\to 0$   as $\delta\to 0^+$.
 This completes the proof of limit~\eqref{1.10}.
  \begin{remark}\label{R:R}~Limit~\eqref{1.9} can be established by  repeating the   arguments of the proof of limit~\eqref{1.10}, by considering the   function
$$
	w(t)= u(\delta t)-q(t),    
$$
 where $u(t)=\RR_{t}(u_0,\de^{-\frac12}\zeta,\de^{-1}\eta)$, $ q(t)=(q_1(t),q_2(t)),$  
\begin{align*}
   q_1(t)&=v_0+t(\eta_1-   B_1(\zeta_1)),  \\
    q_2(t)&=\te_0+t\eta_2,
\end{align*}
$\eta_i = \pi_i \eta$, and $\zeta_1=\pi_1 \zeta$. 
 See Proposition~2.4 in~\cite{nersesyan-2018} for a proof of a  limit similar to~\eqref{1.9} in the case of parabolic equations with polynomially growing nonlinearities.	\end{remark}

\addcontentsline{toc}{section}{Bibliography}
 \bibliographystyle{alpha}

\def\cprime{$'$} \def\cprime{$'$}
  \def\polhk#1{\setbox0=\hbox{#1}{\ooalign{\hidewidth
  \lower1.5ex\hbox{`}\hidewidth\crcr\unhbox0}}}
  \def\polhk#1{\setbox0=\hbox{#1}{\ooalign{\hidewidth
  \lower1.5ex\hbox{`}\hidewidth\crcr\unhbox0}}}
  \def\polhk#1{\setbox0=\hbox{#1}{\ooalign{\hidewidth
  \lower1.5ex\hbox{`}\hidewidth\crcr\unhbox0}}} \def\cprime{$'$}
  \def\polhk#1{\setbox0=\hbox{#1}{\ooalign{\hidewidth
  \lower1.5ex\hbox{`}\hidewidth\crcr\unhbox0}}} \def\cprime{$'$}
  \def\cprime{$'$} \def\cprime{$'$} \def\cprime{$'$}

\end{document}